\documentclass{amsart}
\newif\ifpdf \ifx\pdfoutput\undefined \pdffalse \else \pdfoutput=1 \pdftrue \fi
\ifpdf
    \usepackage[pdftex]{hyperref}
    \usepackage{color,graphicx}
    \DeclareGraphicsExtensions{.pdf,.jpg}
\else
    \usepackage{color,graphicx}  
\fi
\ifpdf
  \newcommand \registercomment[5]
     {\newcommand{#1}[1]{\pdfannot{/Type /Annot /Subtype /Text /T (#3 
) /C #4 /Contents (##1) }}}
\else
  \newcommand \registercomment[5]
     {\newcommand{#1}[1]{\ifhmode\unskip\fi 
\begin{color}{#5}\footnote{\begin{color}{#5}  #2: ##1 
\end{color}}\end{color} \ignorespaces}}
\fi
\renewcommand\registercomment[5]{\newcommand{#1}[1]{\ifhmode\unskip{ }\fi\ignorespaces}}

\registercomment{\HScomment}{HS}{Herman}  {[ 1  1  0 ]}{red}
\registercomment{\BScomment}{BS}{Brigitte}{[ 0  1  1 ]}{magenta}
\registercomment{\FScomment}{FS}{Paco}    {[ 1  0  1 ]}{green}
\registercomment{\DOcomment}{DO}{David}   {[ 0  1  0 ]}{blue}
\registercomment{\GRcomment}{GR}{Guenter} {[.5  1  0 ]}{cyan}
\registercomment{\WWcomment}{WW}{Walter}  {[ 0  1 .5 ]}{black}

\begin{document}
\bibliographystyle{plain}
\newtheorem{theorem}{Theorem}
\newtheorem{lemma}{Lemma}
\newtheorem{conjecture}{Conjecture}
\newtheorem{proposition}{Proposition}
\newtheorem{corollary}{Corollary}
\theoremstyle{definition}
\newtheorem{definition}{Definition}
\theoremstyle{remark}
\newtheorem{open}{Open Problem}

\newcommand{\eps}{\varepsilon}

{\count0=\time \divide\count0 by 60 \count2=\count0
\multiply \count2 by -60 \advance \count2 by \time
\ifnum\count2<10
\xdef\urzeit{\the\count0:0\the\count2}\else
\xdef\urzeit{\the\count0:\the\count2}\fi}
\title{Non-crossing frameworks with non-crossing reciprocals\thanks{}}

\author[Orden]{David Orden$^1$}
\address[David Orden and Francisco Santos]{Departamento de Matem\'aticas,
          Estad\'{\i}stica y Com\-pu\-ta\-ci\'on,
          Universidad de Cantabria,
          E-39005 Santander, Spain.}
\thanks{$^1$ Supported by grant BFM2001-1153, Spanish Ministry of Science 
and Technology}
\email{ordend, fsantos@unican.es}
\urladdr{http://www.matesco.unican.es/\char`\~ordend \char`\~santos}
\author[Rote]{G\"unter Rote$^2$}
\address[G\"unter Rote]{Institut f\"ur Informatik,
          Freie Universit\"at Berlin,
          Taku\-stra{\ss}e 9, D-14195~Berlin, Germany.}
\thanks{$^2$ Partly supported by the Deutsche Forschungsgemeinschaft 
(DFG) under grant RO~2338/2-1.}
\email{rote@inf.fu-berlin.de}
\urladdr{http://www.inf.fu-berlin.de/\char`\~rote}
\author[Santos]{Francisco Santos$^1$}
\author[B. Servatius]{Brigitte Servatius}
\author[H. Servatius]{Herman Servatius}
\address[Brigitte and Herman Servatius]{Dept. of Mathematical Sciences,
          Worcester Polytechnic Institute,
          Worcester, MA 01609}
\email{bservat, hservat@wpi.edu}
\urladdr{http://users.wpi.edu/\char`\~bservat/\ \char`\~hservat/\ }
\author[Whiteley]{Walter Whiteley$^3$}
\address[Walter Whiteley]{Department of Mathematics and Statistics,
          York University, Toronto, Canada.}
\email{whiteley@mathstat.yorku.ca}
\thanks{$^3$ Supported by grants from NSERC (Canada) and NIH (US)}
\urladdr{http://www.math.yorku.ca/Who/Faculty/Whiteley/}


\begin{abstract}
We study non-crossing frameworks in the plane for
which the classical reciprocal on the dual graph is also non-crossing.
We give a complete description of the self-stresses on non-crossing
frameworks $G$ whose reciprocals are non-crossing, in terms of: the
types of faces (only
pseudo-triangles and pseudo-quadrangles are allowed); the sign
patterns in the stress on $G$; and a geometric condition on the stress
vectors at some of the vertices.

As in other recent papers where the interplay of non-crossingness
and rigidity of straight-line plane graphs is studied,
\emph{pseudo-triangulations} show up as objects of special interest.
For example, it
is known that all planar Laman circuits can be embedded as a
pseudo-triangulation with one non-pointed vertex. We show that for
such pseudo-triangulation embeddings of planar
Laman circuits which are sufficiently generic,
  the reciprocal is non-crossing and again a
pseudo-triangulation embedding of a planar Laman circuit. For a
singular (non-generic) pseudo-triangulation embedding of a planar Laman
circuit, the reciprocal is still non-crossing and a
pseudo-triangulation, but its underlying graph may not be a Laman
circuit. Moreover, all the pseudo-triangulations which admit a
non-crossing reciprocal arise as the reciprocals of such, possibly
singular, stresses on pseudo-triangulation Laman circuits.

All self-stresses on a planar graph correspond to liftings to
piece-wise linear surfaces in 3-space.  We  prove
characteristic  geometric properties of the lifts of such non-crossing
reciprocal pairs.

\end{abstract}

\maketitle

\section{Introduction}\label{sectintro}
\subsection{History of Reciprocals}
\begin{figure}[htb]
\centering
\includegraphics[scale=.80]{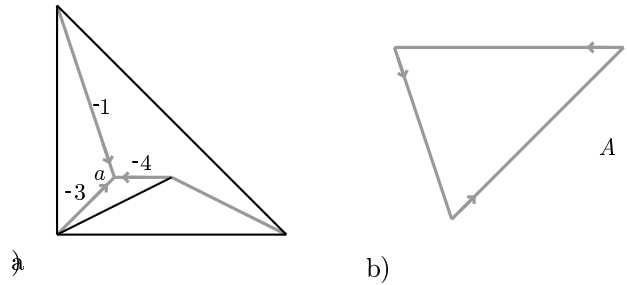}
\caption{The equilibrium forces at vertex $a$ produce a polygon of forces  $A$.
\label{figurewheel4recippseudoonevert}}
\end{figure}

\begin{figure}[htb]
\centering
\includegraphics[scale=.80]{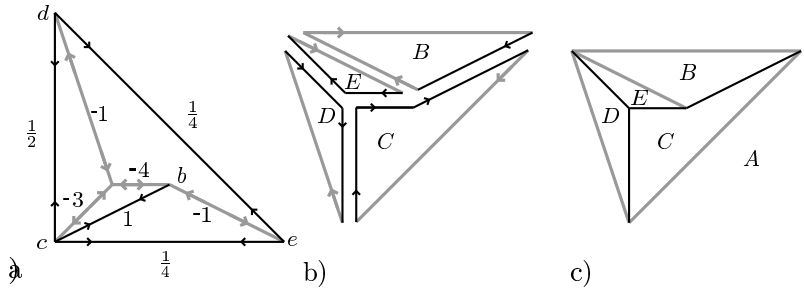}
\caption{Assembling the reciprocal\label{figurewheel4recippseudoeplode}}
\end{figure}

There is a long history of connections between the rigidity of
     frameworks and techniques of drawing planar graphs in the plane.
     Tutte's famous rubber band method~\cite{Tutte3} uses physical
     forces and static equilibrium to obtain a straight line embedding
     of a $3$-connected graph with convex faces, a plane version of
     Steinitz's theorem~\cite{steinitz} that every $3$-connected
     planar graph can be represented as the skeleton of a
     $3$-dimensional polytope.
\begin{figure}[htb]
\centering
\includegraphics[scale=.80]{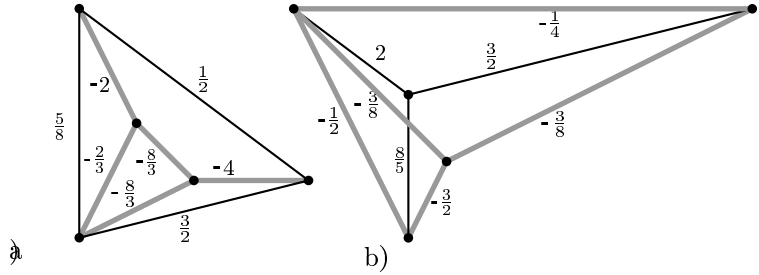}
\caption{A non-crossing framework and crossing reciprocal.
\label{figurewheel4recipnonpseudo}}
\end{figure}
     Maxwell's theory of reciprocal figures~\cite{Maxwell} constructs 
a specific geometric drawing
     of the combinatorial dual of a drawn planar graph provided that 
every edge of the original participates
     in an internal equilibrium stress. However, neither primal nor 
dual are necessarily  crossing free,
     Figure~\ref{figurewheel4recipnonpseudo}.
One can verify the equilibrium of a set of forces at a point by
placing the vectors head-to-tail as a \emph{polygon of forces},
Figure~\ref{figurewheel4recippseudoonevert}.
  For a planar graph with internal forces in equilibrium at each 
vertex, one creates a set of polygons, one
for each vertex.  For each edge of the framework, the two forces at its
endpoints are equal in size and opposite in direction.  With the
polygons sharing parallel edges they can then be pieced together,
edge to edge, as the faces of the dual graph,
Figure~\ref{figurewheel4recippseudoeplode}.
For a framework on a planar graph
with only internal forces (tension and compression in
equilibrium), the creation of a complete reciprocal diagram---a
second framework on the planar dual graph from all these patches---is
equivalent to verifying the equilibrium of these forces.
Each set of equilibrium forces in a planar framework generates
a reciprocal framework, unique up to
translation,
and each reciprocal
prescribes a set of forces in equilibrium.

Reciprocal figures were first developed in the 19th century as a
graphical technique to calculate when external static forces on a
plane framework reach an equilibrium at all the vertices with
resolving tensions and compressions in the
members~\cite{Maxwell,Cremona}, and were the basis for {\em
graphical statics} in civil engineering in the last half of the
19th century. This technique has been rediscovered and applied by
a number of authors, and was used to check the Eiffel tower prior
to construction \cite{Bow,Culmann}.

Classically, there are two graphic forms for the reciprocal.  In
the engineering work adapted to graphical techniques at the
drafting table and presented by the mathematician and engineer
Cremona \cite{Cremona},  the edges of the reciprocal are drawn parallel to
the edges of the original.  We will use this `Cremona' form of the
reciprocal in most of our proofs.

\begin{figure}[htb]
\centering
\includegraphics[scale=.80]{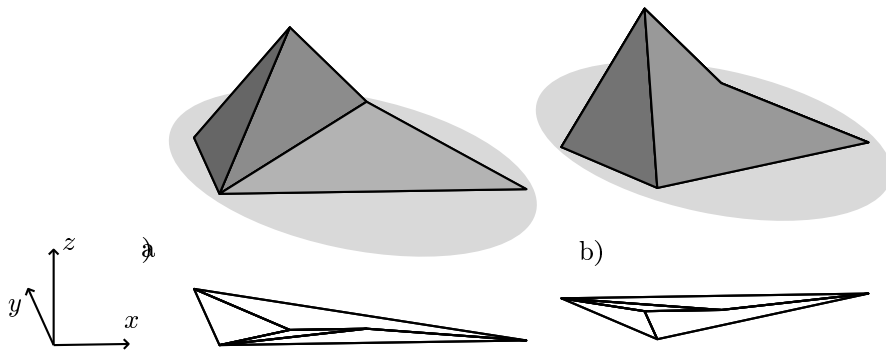}
\caption{The lift of the frameworks of 
Figure~\ref{figurewheel4recippseudoeplode}.\label{figurewheel4lpseudolift}}
\end{figure}

In the original work of Maxwell the edges of the reciprocal
are drawn perpendicular to the edges of the original framework.
This form is adapted to viewing the framework on a planar graph
  as a projection of a spatial polyhedron and the
reciprocal as a drawing of the dual polyhedron.  The surprise is that
this image captures an exact correspondence: a framework with a planar
graph has a self-stress if and only if it is the exact projection
of a spatial, possibly self-intersecting, spherical polyhedron,
Figure~\ref{figurewheel4lpseudolift},
  with non-zero dihedral angles
directly corresponding to the non-zero forces in the self-stress
of the projection \cite{Maxwell,Cremona,Whiteley2,cw01,WhiteWhiteley2}.
Moreover, the reciprocal diagram is the projection of
a very specific spatial polar of the original projected
polyhedron.  We will return to an application of this
correspondence with spatial polyhedra in Section~\ref{sectlift}.  This spatial polarity also reinforces the 
reciprocal relationship for
plane frameworks, since either framework can be viewed as the original and
the other viewed as its reciprocal.
It is not difficult to check that given
one presentation of the reciprocal, we can simply turn it $90^\circ$
to create the other presentation.

As a modern connection, we note that reciprocal diagrams were
rediscovered as a technique to check whether a given plane drawing
is the exact projection of a spatial polyhedron or a polyhedral
surface~\cite{Huffman, Huffman2}.  Here the reciprocal diagram is also
called the `gradient diagram', since the vertices of the
reciprocal can be located as the points of intersection  of normals to the
faces of the original polyhedron, with the
projection plane, with all normals drawn from a fixed center above the
plane. These points are the gradients representing the slopes of the faces.
This gradient diagram is a Maxwell reciprocal with
reciprocal edges perpendicular to the edges of the original
\cite{cw01}. A related construction, starting with points on the
paraboloid $x^2 + y^2 -z = 0$ and their convex hull, creates the
Delaunay triangulation of the projected points, with the Voronoi cells as a
(Maxwell) reciprocal diagram \cite{abcw}.

\subsection{Our Contribution}
Here we pull together planar embeddings with the theory of
reciprocals: we investigate when both the original graph
drawing, or framework, and the reciprocal diagram are crossing
free.  We show in Section~\ref{sectsimultaneously}
that interior faces of both drawings must be
pseudo-triangles or pseudo-quadrangles, the outer boundary must be convex, 
and the self-stress generating the reciprocal must have a specific type
of sign pattern.  With one additional geometric condition on the
self-stress at non-pointed vertices, these conditions become
both necessary and sufficient for creating a reciprocal pair
of non-crossing frameworks. The cases where the framework is a Laman
circuit (an edge-minimal graph that can sustain a self-stress in a generic
embedding), a pseudo-triangulation, or has a unique non-pointed vertex,
are specially interesting, and addressed separately.

There is a case that combines the three just mentioned: that of Laman
circuits embedded as pseudo-triangulations. Previous work
\cite{horsssssw} shows that all planar Laman circuits admit such
embeddings, and they have a unique non-pointed vertex. In
Section~\ref{sectsimultcirc} we show that they produce
non-crossing reciprocals.  For sufficiently generic
embeddings the reciprocals are also Laman circuits realized as
pseudo-triangulations with one non-pointed vertex, so the pairing is a
complete reciprocity. This result is extended to the Laman circuits
realized as non-generic pseudo-triangulations, with self-stresses that
are zero on some edges, again showing that the reciprocal is a
pseudo-triangulation (with several non-pointed vertices, in general).
Moreover, any pseudo-triangulation with a self-stress and a
non-crossing reciprocal occurs in a pair with a possibly singular
pseudo-triangulation on a Laman circuit.

Section~\ref{sectlift} studies the characteristic features of the 
lifts of such non-crossing reciprocal pairs. Essentially, they look
like negatively curved surfaces with a unique singularity at the
vertex whose reciprocal is the outer face. In particular, that vertex
is the unique local maximum of the surface, the outer face is the unique
local minimum, and there are no (horizontal) saddle-points.

We finish with a list of open problems on questions related to these, 
in Section~\ref{sectopen}.


\subsection{Preliminaries---Frameworks}
     An \emph{$n$-dimensional framework} $(G,\rho)$ is a graph $G$
together with an embedding $\rho\colon V \rightarrow \mathbb{R}^n$,
and we will write $\rho(i) = \mathbf{p}_i$.  In this paper we will
only consider frameworks in the plane.  The edges of $(G,\rho)$ are
regarded as abstract length constraints on the motions of
$\mathbf{p}_i$.  The edges are often drawn in as straight line
segments, which may, of course, cross.  In the absence of any edge
crossing, we will say that the framework is \emph{non-crossing}.

     Infinitesimally, distance constraints form a linear system, with 
an equation
for each edge $\{i,j\}\in E$
     $$(\mathbf{p}'_i-\mathbf{p}'_j)\cdot(\mathbf{p}_i-\mathbf{p}_j) =
{0}$$
  This system of equations in the unknowns $\mathbf{p}'_i$ has
a coefficient matrix of size  $|E|\times 2|V|$, the \emph{rigidity matrix} of
the framework.
The framework  is called \emph{infinitesimally
rigid} if its
  rigidity matrix has rank $2|V|-3$.
The case where $G$ is the complete graph deserves special attention.
Then the rigidity matrix
has size $\binom{|V|} {2}\times 2|V|$ and rank $2|V|-3$
(unless all the vertices lie on a single line).
The matroid of (rows of) the rigidity matrix of the complete graph,
called
the \emph{rigidity matroid} \cite{gss,Whiteley3} of the point set
$\rho(V)$, is interesting because its spanning subsets are precisely
the infinitesimally rigid frameworks with vertex set $\rho(V)$.

The rigidity matroid is the same for all generic choices of vertex positions,
and is called the \emph{generic rigidity matroid}. Spanning graphs of it
are called \emph{generically rigid} graphs and the minimal ones (bases
of the matroid) are called \emph{isostatic} or \emph{Laman graphs}.
They are characterized
by the \emph{Laman condition}: $G=(V,E)$ is isostatic if and only if
$|E| = 2|V|-3$ and every subset of $k\geq 2$
     vertices spans at most $2k-3$ edges of $E$, \cite{Laman}.  Generically
rigid graphs are those containing a spanning  \emph{Laman subgraph}.
     Every circuit of the generic rigidity matroid which spans $k$ vertices
must consist of exactly
     $2k-2$ edges, and is called a \emph{Laman circuit.}

     In a dual analysis of this matrix, a stress on a framework $(G,\rho)$ is an
assignment
     of scalars $\omega\colon E\rightarrow \mathbb{R}$.   A stress is 
\emph{resolvable}
or a \emph{self-stress} of the framework
     if the weighted sum of the displacement vectors corresponding to 
each vertex
cocycle is zero;
     $$\sum_{j|(i,j)\in E}  \omega_{ij} (\mathbf{p}_i 
-\mathbf{p}_j)=\mathbf{0} ,
    \mbox{~for all ~$i \in V$.}$$
That is to say, the self-stresses form the cokernel of the rigidity matrix.
  If the graph is generically rigid and the embedding is generic,
  then the dimension of the space of self-stresses is
  $|E|-(2|V|-3)$.  In particular, a Laman circuit in a generic embedding has
a unique (up to a scalar multiple) self-stress, which is
non-zero on
     every edge.

\subsection{Preliminaries---Reciprocal diagrams}
A \emph{plane graph} $G\to\mathbb{R}^2$ is a graph which is
(topologically) embedded in the plane.  The embedding determines the
combinatorial information about the sequences of edges that lie on the
boundary of each face (the \emph{face cycles}) and the sequences of
edges that lie around each vertex (the \emph{vertex cycles}).  A plane
graph $G\to\mathbb{R}^2$ determines a \emph{dual plane graph}
$G^*\to\mathbb{R}^2$ which has a vertex for each face of
$G\to\mathbb{R}^2$, an edge between two vertices for each edge
separating the corresponding faces of $G\to\mathbb{R}^2$, and a
face for every vertex of $G$.  Vertex cycles of $G\to\mathbb{R}^2$
correspond to face cycles of $G^*\to\mathbb{R}^2$ and vice
versa. The dual graph is unique up to choice of which vertex of $G$
will become the unbounded face of $G^*\to\mathbb{R}^2$.

Given a plane graph $G\to \mathbb{R}^2$ and a framework $(G,\rho)$
on $G$, a second framework on the plane dual graph $G^*$ is
\emph{reciprocal} to the first if corresponding edges are parallel.
Even if the framework $(G,\rho)$ is already non-crossing, one may choose
to use a different plane embedding of $G$ for computing the reciprocal
(of course, unless $G$ is 3-connected in which case its plane
embedding is unique up to choice of the outer face and orientation). But in
this paper we will only consider the case where $(G,\rho)$ is non-crossing and
$G\to\mathbb{R}^2$ is the embedding given by $\rho$. In particular, we
omit  mention of the plane embeddings $G\to \mathbb{R}^2$ and
$G^*\to \mathbb{R}^2$ in the sequel, since both can be deduced
from the (non-crossing) framework $(G,\rho)$.
We say that a reciprocal $(G^*,\rho^*)$ is \emph{a non-crossing
reciprocal} of $(G,\rho)$ if $(G^*,\rho^*)$ is non-crossing \emph{and}
its embedding is dual to the embedding of $(G,\rho)$. Our goal is to
characterize pairs of simultaneously non-crossing reciprocal
diagrams. As a first (counter-)example,
Figure~\ref{figurewheel4recipnonpseudo} shows a non-crossing framework
with crossing reciprocal, which is actually the ``typical
situation''.

Observe that, in principle, if a non-crossing graph has a non-crossing
reciprocal, the reciprocity may preserve or reverse the
orientation. That is, vertex cycles of $(G,\rho)$ may in principle become face
cycles in $(G^*,\rho^*)$ in the same or the opposite directions. We
will prove that only the orientation-reversing situation occurs.

\medskip

Reciprocity of frameworks is very closely related to self-stresses.
We offer a simple representation of the Cremona reciprocal for a plane
graph~$G$.
We denote the set of all directed edges of $G$
and their inverses by $E^\pm$. The inverse of an edge $e$ is denoted
by $\bar e$, with $\bar{\bar e}=e$.  Let $g \colon E^\pm \rightarrow
\mathbb{R}^2$ be an assignment of unit vectors to the edges of $G$
with $g({\bar e}) = -g(e)$, $|g(e)| = 1$.  For simplicity we will
write $g(e) = \mathbf{e}$.
%
Suppose we have two scalar functions $\alpha\colon E \rightarrow
\mathbb{R}$ and $\beta\colon E \rightarrow \mathbb{R}$ with
compatibility conditions,
\begin{equation} \label{cyclecocyclecondition}
       \sum_{e \in C} \alpha_e \mathbf{e} = \mathbf{0}, \qquad
       \sum_{e \in C'} \beta_e \mathbf{e} = \mathbf{0}
\end{equation}
for each facial cycle $C$ and each vertex cocycle $C'$.

Since the facial cycles 
corresponding to a plane
graph embedding generate the entire cycle space of the graph, the cycle
conditions in~(\ref{cyclecocyclecondition}) are sufficient to
guarantee that the displacement vectors $\alpha_e \mathbf{e}$ are
consistent over the entire framework. Hence, they are the
displacements or edge vectors of a framework $(G,\rho)$ on the graph $G$.
Similarly,
the vectors $\beta_e \mathbf{e}$ correspond to edge displacements of
a framework $(G^*,\rho^*)$ on the dual graph $G^*$.
In particular, the two frameworks are reciprocal
to each other. But now, the face equalities for $G^*$ can be read
as equilibrium conditions for the vertices of $(G,\rho)$ and
vice versa. Hence,
if $\alpha_e \neq 0$ then the values
$\beta_e/\alpha_e$ are a self-stress on the framework~$(G,\rho)$.
Similarly, if $\beta_e \neq 0$ then the values
$\alpha_e/\beta_e$ are a self-stress on the framework~$(G^*,\rho^*)$.

This argumentation can be reversed, and a reciprocal framework can be
constructed uniquely starting with a self-stress on $(G,\rho)$ 
\cite{cw01,cw02}.
It follows that the self-stresses of a connected framework $G$ are in
one-to-one correspondence with the reciprocals of a given framework
$G$ (up to translation). Multiplication of the self-stress by a
constant corresponds to scaling the reciprocal. In particular,
changing the sign of a self-stress will rotate the reciprocal by
$180^\circ$. Thus, if the framework $G$ has a unique self-stress (up
to scalar multiplication), we can speak of \emph{the} reciprocal
framework if we are not interested in the scale.

Maxwell proved that the projection of a spherical polyhedron from
3-space gives a plane diagram of segments and points which forms a
stressed bar and joint framework. This proof, and related constructions, were
built upon an analysis of reciprocal diagrams \cite{Maxwell}.  Crapo and
Whiteley~\cite{cw02,
   Whiteley2} gave new proofs for Maxwell's theorem as well as the
converse for planar graphs.
See~\cite{cw01} and~\cite{cw02} for more details on the full vector spaces
of self-stresses, reciprocals and spatial liftings of a plane
drawing.


\subsection{Preliminaries---Pseudo-triangulations}
\begin{figure}[htb]
\centering
\includegraphics[scale=.80]{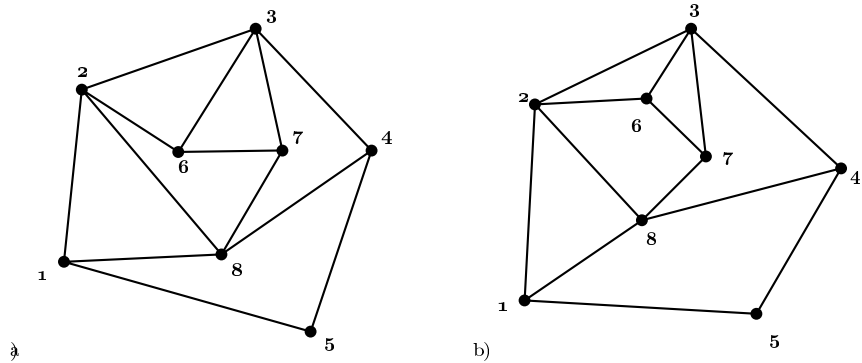}
\caption{(a) A pointed pseudo-triangulation (b) a different embedding
of the same graph which is not a pseudo-triangulation.\label{figurepseudot}}
\end{figure}

Given a non-crossing embedding $\rho\colon V \rightarrow \mathbb{R}^2$
  of a planar graph, we say that the vertex $i$ is \emph{pointed} if all
adjacent points $\mathbf{p}_j$
lie strictly on
one side of some line through $\mathbf{p}_i$. In this case some pair 
of consecutive edges
in the counter-clockwise order around $i$ spans a \emph{reflex
angle}.  A face of the non-crossing framework is a \emph{pseudo-triangle}
if it is a simple planar polygon with exactly three convex vertices
(called \emph{corners}). A
\emph{pseudo-triangulation} has all interior faces pseudo-triangles and the
complement of the outer face is a convex polygon. See
Figure~\ref{figurepseudot}.
In a \emph{pointed pseudo-triangulation} all the vertices are pointed.

In this paper we need to extend the concept of pseudo-triangulation to that of
\emph{pseudo-quadrangulation}. A \emph{pseudo-quadrangle} is a
simple polygon with four convex vertices (corners) and a pseudo-quadrangulation
is a decomposition of a convex polygon into pseudo-triangles and 
pseudo-quadrangles.

\begin{lemma}
\label{lemmaquadrangulation}
Let $T$ be a pseudo-quadrangulation with $e$ edges,
  $x$ non-pointed vertices,
  $y$ pointed vertices,
  $t$ pseudo-triangles,
  and $q$  pseudo-quadrangles. Then,
\[
2e = t + 3y+4x-4.
\]
\end{lemma}

\begin{proof}
The pseudo-quadrangulation has $3t+4q$ convex angles, and $y$ reflex
angles (one at each pointed vertex). Since the total number of angles
is $2e$, we get $2e=3t+4q+y$. Euler's formula
gives $e=x+y+t+q-1$. Eliminating $q$ gives the desired equation.
\end{proof}

In the case of pseudo-triangulations ($q=0$) this lemma is well-known
and usually stated under the equivalent (via Euler's formula, $t=e+1-x-y$)
form $e=(2n-3) + x$, where $n=x+y$ is the total number of vertices 
\cite{horsssssw,osss}.
Since every pseudo-triangulation is infinitesimally rigid \cite{os}, 
this formula says that:

\begin{lemma}
\label{lemmadimensionpt}
The dimension of the space of self-stresses of a pseudo-triangulation 
equals its number
of non-pointed vertices, that is, $e-(2n-3)$.
\end{lemma}

Moreover, it is easy to prove that every non-crossing framework
can be extended to a pseudo-triangulation with exactly the same 
number of non-pointed vertices~\cite[Theorem~6]{rwwx}. The next lemma follows.

\begin{lemma}
\label{lemmadimension}
The dimension of the space of self-stresses of a non-crossing 
framework is at most its number
of non-pointed vertices.
\end{lemma}

Pseudo-triangulations have arisen as important objects connecting rigidity
and planarity of geometric graphs. For example,
pointed pseudo-triangulations were an important tool in straightening 
the carpenter's
rule~\cite{streinu}.
A graph is planar and generically rigid if and only if
it can be embedded as a pseudo-triangulation~\cite{osss}. It is
planar and isostatic if and only if it can be embedded as a
pointed pseudo-triangulation~\cite{horsssssw}.

These embedding results have extensions for Laman circuits.   A {\it
pseudo-trian\-gu\-lation circuit} is a planar Laman circuit embedded as a
pseudo-triangulation.  A pseudo-triangulation circuit has a single non-pointed vertex,
by Lemma~\ref{lemmadimensionpt}.

\begin{figure}[htb]
\centering
\includegraphics[scale=.85]{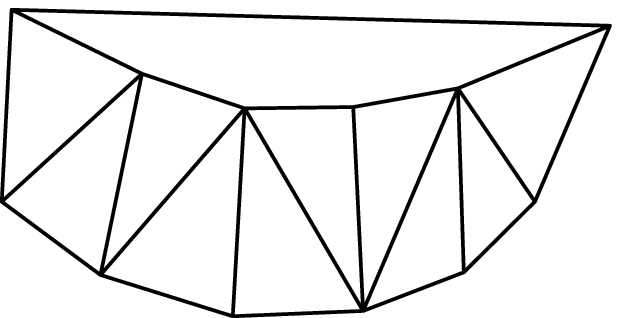}
\includegraphics[scale=.85]{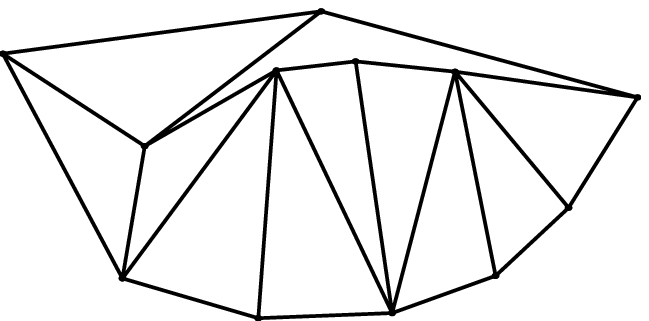}
\caption{A Laman circuit and one of its embeddings as a
pseudo-triangulation with exactly one non-pointed vertex.\label{figurecircuit}}
\end{figure}

See Figure~\ref{figurecircuit} for an example of a Laman circuit (a Hamiltonian
polygon triangulation with an added edge between its two vertices of 
degree $2$)
and one of its embeddings as a pseudo-triangulation with exactly one 
non-pointed
vertex. A basic starting point for our analysis is the following result from
\cite{horsssssw}.

\begin{theorem} \cite{horsssssw} Every topologically embedded planar
   Laman circuit with a given outer face and a specified vertex $v$
   that does not lie on the outer face has a realization as a
   pseudo-triangulation where $v$ is the single non-pointed vertex.
\label{theorem2connected}
\end{theorem}

These pseudo-triangulation circuits are the main focus of
Section~\ref{sectsimultcirc}.

\subsection{Geometric versus Singular Circuits}
\begin{figure}[htb]
\centering
\includegraphics[scale=0.675]{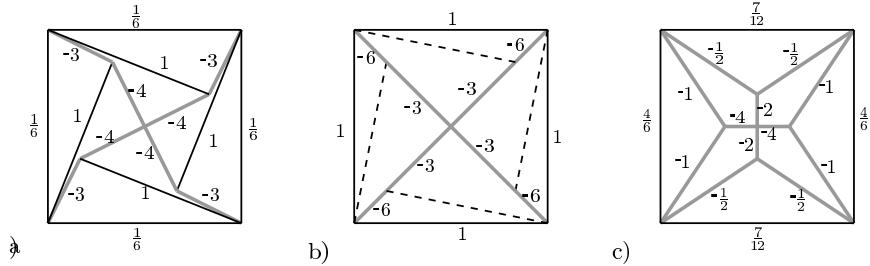}
\caption{A geometric circuit (a), a singular realization with dropped edges
(b) and a singular realization with additional
self-stresses (c).
\label{figuresingcyc}}
\end{figure}

Given a planar Laman circuit $G$,
we have a range of realizations as frameworks in the
plane, all of which are dependent, i.e. have  a non-trivial space of
self-stresses.
   For an open dense subset of these realizations, containing the generic
realizations, the unique (up to scalar multiplication) self-stress
  is non-zero on all edges.  We say the graph is embedded as a {\em
geometric circuit}, see Figure~\ref{figuresingcyc}(a).

The remaining {\em singular realizations\/} are frameworks on which either the
one-dimensional space of self-stresses vanishes on some subset of edges,
as in Figure~\ref{figuresingcyc}(b) in which the dashed edges are unstressed,
or for which the space of self-stresses has higher dimension,  in which case
it is generated by self-stresses which vanish on some of the edges.  See
Figure~\ref{figuresingcyc}(c)
self-stress. The singular self-stresses of Figure~\ref{figuresingcyc}(b)
  will cause some complications in
the reciprocal diagrams: the original edge effectively disappears as a
division between faces and the corresponding reciprocal edge has zero
length, fusing the reciprocal pair of vertices into one.  We can actually
track such singular frameworks $\rho$ on the graph as those whose vertex
coordinates satisfy at least one of a set of $e$ polynomials,
$C_{i,j} (\rho)$, representing the pure conditions for the independence
of the sub-graphs with the edge $i,j$ removed ~\cite{WhiteWhiteley2}.  In
general, the coefficients of the unique self-stress of a geometric circuit on
the original graph can be written using these polynomials as coefficients.  An
edge has a zero coefficient in the self-stress if and only if  the 
corresponding
polynomial is zero.

\begin{figure}[htb]
\centering
\includegraphics[width=2.0in]{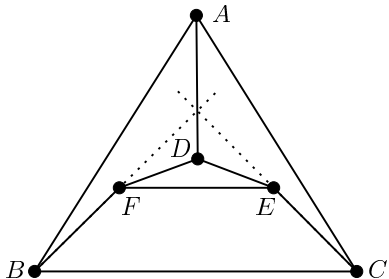}
\caption{A Laman circuit with a singular self-stress.}
\label{figurepseudoNonCircuit}
\end{figure}

However, the realizations as pseudo-triangulations are not guaranteed to be
geometrical circuits, since they need not be generic embeddings.
In the pseudo-triangulation of Figure~\ref{figurepseudoNonCircuit},
the  edge $CF$
does not participate in the self-stress when the edges $AD$, $BF$, and $CE$
are concurrent, for projective geometric reasons.
By Lemma~\ref{lemmadimensionpt}, all pseudo-triangulation realizations
of this graph will have a $1$-dimensional space of self-stresses,
but the stress may be singular.


\section{Simultaneously non-crossing reciprocals}\label{sectsimultaneously}

Assume we are given a non-crossing  framework $(G,\rho)$, and a 
particular everywhere
non-zero self-stress to construct the reciprocal from.
The goal of this section is to determine the conditions for the
framework and its reciprocal to be simultaneously non-crossing.
  Our main result is that for this to happen the framework needs to
be a pseudo-quadrangulation
and there are certain necessary (and almost sufficient)
conditions on the
signs of the self-stress that make the reciprocal non-crossing
(Theorem~\ref{theoremsimultaneous}).

\begin{figure}[htb]
\centering
\includegraphics{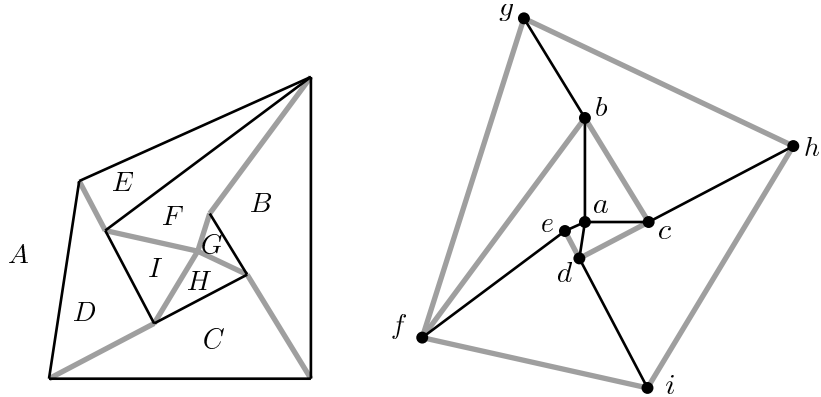}
\caption{We can have flat angles on both sides of a non-crossing
reciprocal pair.}
\label{figureflatreciprocal}
\end{figure}

\begin{figure}[htb]
\centering
\includegraphics[scale=.80]{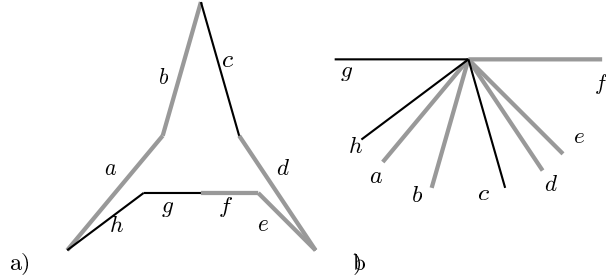}
\caption{Locally, the reciprocal of a
pseudo-quandrangle with a flat vertex as a corner is a flat vertex. }
\label{figureflatlocal}
\end{figure}

In order to include all degeneracies that arise in
non-crossing reciprocal pairs, we do not assume our framework to be in general
position. In particular, angles of exactly $\pi$ (to be called 
\emph{flat angles}) can arise and, by convention,
we treat them as convex (``small'') angles. In particular, a vertex 
having one such flat angle
is necessarily non-pointed, and the face incident to it cannot be a 
pseudo-triangle.
See Figures~\ref{figureflatreciprocal},~\ref{figureflatlocal}.
  We do not allow degeneracies in which the angle is
$0$, as these produce `crossing' edges. Also, vertices with two flat 
angles have degree two and
produce a double edge in the reciprocal, so we do not allow them in 
our frameworks.

Recall that for the reciprocal of $(G,\rho)$ to be considered non-crossing,
the face cycles in  $(G,\rho)$ must form the vertex cycles in $(G^*,\rho^*)$.
In principle, the face cycles
could become vertex cycles with either the same or the opposite orientation,
but the orientation must be globally consistent: either all the cycles keep
the orientation or all of them reverse it.
Special care needs to be taken with the exterior
face, for which the orientation must be considered ``from outside'':
if the orientation
is chosen counterclockwise for interior faces, it will
be clockwise for the outer face.

\subsection{The reciprocal of a single face}
\label{subsecsingleface}

Let us concentrate on a single face of our framework $G$. That is
to say, let $F$ be a simple polygon in the plane. If $F$ has $k$
convex vertices we say it is a pseudo-$k$-gon.

Given a vertex $v$ of $F$, we call the \emph{reduced internal angle of $F$ at
$v$} the internal angle itself if $v$ is a convex vertex (corner) of $F$,
and the angle minus $\pi$ if it is not.
In other words, we are
``reducing'' all angles to lie in the range $0 < \theta \leq \pi$.
Our first result
generalizes the elementary fact that the total
internal angle of a $k$-gon is $(k-2)\pi$.


\begin{lemma} \label{lemmatotalredintangle} The sum of reduced internal
angles of a pseudo-$k$-gon is $(k-2)\pi$.
\end{lemma}

\begin{proof}
If the polygon has $n$ vertices,
the sum of (standard) internal angles is $(n-2)\pi$, and the
reduction process subtracts $(n-k)\pi$.
\end{proof}


Now let signs be given to the edges of $F$, intended to represent
the signs of a self-stress in the framework of which $F$ is a face.
In this section we assume no sign is zero, which is not really a loss
of generality: if a self-stress is zero on some edges then it is a
self-stress on the subframework on which it is not zero and the reciprocal
depends only on that subframework.

In the reciprocal framework, $F$ corresponds to a certain vertex to which
the reciprocal edges are incident.
In this section we will use the following exact rule to draw the
 reciprocal edges:\GRcomment{Most figure (even in this section) use the opposite rule!
We might consider to switch this one (in the next version).}
walk along the boundary of $F$ such that $F$ is on the left side of the edges,
i.~e., surround $F$ counter-clockwise.
Then, edges with positive stress will produce reciprocal edges pointing
in the same direction as the original, and edges with negative stress
will produce reciprocal edges pointing in the opposite direction at that
vertex.
These conventions are no loss of generality;
the opposite choice would produce a reciprocal rotated by 180 degrees.
  The actual value of the stress will give the length of the
reciprocal edges.

What are the conditions for the reciprocal to be ``locally'' non-crossing?
The reciprocal edges must appear around the
reciprocal vertex of $F$ in the same (if we want an orientation-preserving
reciprocal) or opposite (for an orientation-reversing reciprocal)
cyclic order as they appear in $F$. That is to say, if $e_1$ and $e_2$ are
two consecutive
edges (in counter-clockwise order) of $F$ and $e^*_1$ and $e^*_2$ are the
reciprocal edges, we want the angle from $e^*_1$ to $e^*_2$ in the reciprocal
to contain no other edge, where the reciprocal
angle is taken counter-clockwise
if we want to preserve orientation and clockwise if we want to reverse it.
A necessary and sufficient condition for this to happen is that
\emph{the
sum of angles
\textup(measured all clockwise or all counter-clockwise,
depending on the orientation case we are in\textup)
between reciprocals of consecutive edges of $F$ add up to $2\pi$}.
We are now going to translate this into a condition on the signs of edges.

Let us first look at the orientation-reversing case.
Let $e_1$ and $e_2$ be two consecutive edges of $F$ with
common vertex~$v$. We say that the angle at $v$ has \emph{face-proper}
signs 
(or that the angle of $F$ at $v$ is face-proper,
for short)
if either
$v$ is a corner of $F$ and the signs of $e_1$ and $e_2$
are opposite, or $v$ is not a corner and the two signs are equal.
When the signature is not face-proper (a corner with no sign change or a
reflex-angle with sign change) we call it \emph{vertex-proper}.
The reason for this terminology is that
when rotating a ray around a vertex, the fastest (``proper'')
way of going from an edge to the next one is to keep the direction of the ray
for a convex angle and to change to the opposite ray (``changing signs'')
when the angle is reflex.
Analogously, when sliding a tangent ray around a polygon, we should
change to the opposite direction at corners.

The key fact now is that  the reciprocal angle
(measured clockwise) of a given vertex $v$ of $F$ equals the reduced 
internal angle at
$v$ if the signature is face-proper, and it equals the reduced internal angle
plus $\pi$ if it is vertex-proper.
Hence:

\begin{lemma}
\label{lemmareverse}
The sum of angles \textup(measured all clockwise\textup)
between reciprocals of consecutive edges of $F$ equals the total
reduced internal angles of $F$ plus $\pi$ times the number of
vertex-proper angles.

In particular, in order for $F$ to produce a planar reciprocal
with the orientation reversed, $F$ must be either a pseudo-quadrangle
with no vertex-proper angle, or a pseudo-triangle
with only one vertex-proper angle.
If there is a flat angle $\pi$,
it must occur in a pseudo-quadrangle with a sign change at
all corners \textup(one of which is the flat angle\textup).
\end{lemma}

\begin{proof}
The two cases described are the only ways of getting $(k-2)\pi +s\pi=2\pi$,
where $k$ is the number of corners (and hence $(k-2)\pi$ is the reduced
internal angle of $F$) and $s$ the number of vertex-proper angles.

If one of the convex angles is flat, then there must be four convex
angles to achieve a non-crossing polygon.  This is a pseudo-quadrangle
with a sign change at each of the corners.
\end{proof}

\begin{figure}[htb]
\centering
\includegraphics{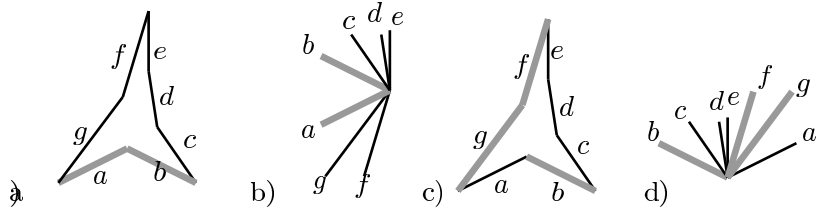}
\caption{(a) a pseudo-triangle with one vertex-proper angle, at the 
corner between $e$ and~$f$ and (b) its reciprocal vertex,
(c) a pseudo-triangle with one vertex-proper angle between $a$ and $b$,  not at a
corner, and (d) its reciprocal vertex.\label{figureorden0822bc}}
\end{figure}

\begin{figure}[htb]
\centering
\includegraphics{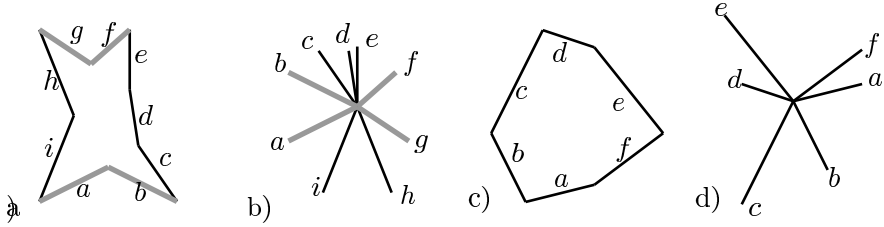}
\caption{A pseudo-quadrangle with all angles face-proper (a) and
a convex face with no sign changes (c), with dual 
vertices\label{figureorden0822da}}
\end{figure}

Figure~\ref{figureorden0822bc}
and the left half of Figure~\ref{figureorden0822da} illustrate the 
cases permitted by
Lemma~\ref{lemmareverse}. Thick and thin lines represent the
two different signs. The case of a pseudo-triangle produces two different
pictures, depending on whether the vertex-proper angle happens at a corner
(Figure~\ref{figureorden0822bc}a)
of the pseudo-triangle or at a reflex vertex
(Figure~\ref{figureorden0822bc}c).
Parts (b) and (d) show the reciprocal vertex with its
incident edges.

The conditions for the orientation-preserving case are now easy to derive,
and illustrated in Figure~\ref{figureorden0822da}c and 
~\ref{figureorden0822da}d.

\begin{lemma}
\label{lemmapreserve}
Let $n$ be the number of vertices in $F$.
The sum of angles (measured all counter-clockwise)
between reciprocals of consecutive edges of $F$ equals $2n\pi$ minus
the total
reduced internal angle of $F$ and minus $\pi$ times the number of vertex-proper
sign changes.

In particular, in order for $F$ to produce a planar reciprocal
with the same orientation, $F$ must be a strictly convex polygon and all
edges must have the same sign.
\end{lemma}

\begin{proof}
The first assertion follows from Lemma~\ref{lemmareverse} and
the fact that the clockwise and
counter-clockwise angles between two edges add up to $2\pi$.

For the second assertion, the equation that we now have is
$2n\pi - (k-2)\pi - s\pi=2\pi$, where $k$ and $s$ are again the number
of corners and vertex-proper sign changes. This equation reduces to
$s+k=2n$, which implies $s=k=n$. This corresponds to a convex polygon
with all signs equal, as stated.

In particular, an angle of $\pi$ with no sign change will
make the two reciprocal edges overlap at the reciprocal vertex.
Any angle of
$\pi$ with a non-crossing reciprocal will have a sign change.
This means such angles cannot occur in the reciprocal with the
same orientation.   The polygon is strictly convex.
\end{proof}

\subsection{Combinatorial conditions for a non-crossing reciprocal}

Observe that in the description above, $F$ was implicitly assumed
to be an interior face. As mentioned before, orientations of the outer
face have to be considered reversed, which means that the conditions
of Lemmas \ref{lemmareverse} and \ref{lemmapreserve} have to be
interchanged when looking at the outer face. Hence:

\begin{theorem}{\label{nopreservingreciprocals}}
It is impossible for a non-crossing framework to have a
non-crossing reciprocal with the same orientation.
\end{theorem}

\begin{proof}
According to Lemma \ref{lemmapreserve}, all interior faces should be
convex and all edges should have the same sign. But, according to
Lemma \ref{lemmareverse}, the outer face would need to have
some sign changes: at the (at least three) convex hull vertices,
keeping the sign produces a vertex-proper sign change, and we are only
allowed to have one of them.
\end{proof}

Hence, every pair of non-crossing reciprocals will have the
orientations reversed one to the other. In this case,
we have the following statement that follows directly from
Lemmas~\ref{lemmareverse} and~\ref{lemmapreserve}:

\begin{theorem}[\bf Face conditions for a planar reciprocal]
\label{theorem-face-conditions}
    Let $(G,\rho)$ be a non-crossing framework with given self-stress $\omega$.
    The following \emph{face conditions} on the signs of $\omega$ are 
necessary in order for
the reciprocal framework
    $(G^*,\rho^*)$ to be also non-crossing:
\begin{enumerate}
\item the \textup(complement of\textup) the exterior face is strictly 
convex with no sign
changes.
\item the internal faces of $(G,\rho)$ are either
\begin{enumerate}
\item pseudo-triangles with two sign changes, both occurring at corners.
\item pseudo-triangles with four sign changes, three occurring at corners.
\item pseudo-quadrangles with four sign changes, all occurring at corners.
\qed
\end{enumerate}
\end{enumerate}
\end{theorem}

Theorem~\ref{theorem-face-conditions}  says in particular that if two
reciprocal frameworks are both non-crossing, then they are both
pseudo-quadrangulations.

Of course, the conditions on faces of $(G,\rho)$ translate into
conditions on the vertices of $(G^*,\rho^*)$. For both
to be non-crossing, both sets of conditions must be satisfied in
both. In particular, the following vertex conditions be satisfied in 
$(G,\rho)$:

\begin{theorem}[\bf Vertex conditions for a planar reciprocal]
\label{theoremsimultaneous}
\label{theorem-vertex-conditions}
    Let $(G,\rho)$ be a non-crossing framework with given self-stress $\omega$.
    Then, in order for
the reciprocal framework $(G^*,\rho^*)$ to be also non-crossing,
the following \emph{vertex conditions} need to
be satisfied by the signs on its vertex cycles:
\begin{enumerate}
\item there is a non-pointed vertex with no sign changes.
\item all other vertices are in one of the following three cases:
\begin{enumerate}
\item pointed vertices with two sign changes, none of them at the big angle.
\item pointed vertices with four sign changes, one of them at the big angle.
\item non-pointed vertices, including any vertices with a flat angle, with
four sign changes.
\end{enumerate}
\end{enumerate}
Moreover, vertices of $(G,\rho)$ in each of the cases \textup{(1)}, 
\textup{(2.a)},
\textup{(2.b)} and \textup{(2.c)} correspond respectively to faces of
$(G^*,\rho^*)$ in the same parts of
Theorem~\ref{theorem-face-conditions}, and vice versa.
\end{theorem}

\begin{proof}
Each of the cases of Theorem~\ref{theorem-face-conditions}, applied
to a face in $(G^*,\rho^*)$, gives the
condition stated here for the reciprocal vertex of $(G,\rho)$.
See Figures~\ref{figureorden0822bc}--\ref{figureorden0822da}.
\end{proof}

We will show below that the vertex conditions of
Theorem~\ref{theorem-vertex-conditions} actually imply the face
conditions of Theorem~\ref{theorem-face-conditions}.

Both the face and the vertex conditions admit a simple rephrasing in
terms of vertex-proper and face-proper angles. Namely:
\FScomment{Added this}
\begin{enumerate}
\item The face conditions are that there is exactly one vertex-proper
angle in every pseudo-triangle, and no vertex-proper angle in the
pseudo-quadrangles and the outer face.
\item The vertex conditions are that there are exactly three face-proper
angles at every pointed vertex and four at every non-pointed vertex
other than the one reciprocal to the outer face, which has no
face-proper angle.
\end{enumerate}

Observe that a pseudo-$k$-gon with $k$ even (resp., $k$ odd) must have an
even (resp., odd) number of vertex-proper angles, simply because it 
has an even number of sign changes.
Hence, the face conditions are that 
this number is as small as possible:
0 for pseudo-quadrangles and the outer face, 
1 for pseudo-triangles.

Similarly, the number of face-proper angles around a non-pointed 
(resp., pointed) vertex must be even (resp., odd).
Again, the vertex conditions are that this number is as small as
possible for the distinguished non-pointed vertex and for 
all pointed vertices (but not for other non-pointed vertices), 
as the following result shows:

\begin{lemma}
\label{lemmapointed}
A self-stress produces at least
three face-proper angles at every pointed vertex $v$,
\textup(unless it is zero on all edges incident to~$v$\textup).
\end{lemma}

\begin{proof}
Observe that in order to meet the equilibrium condition around a vertex,
no line can separate 
the positive edges of the self-stress incident to that vertex
from the negative ones. In the case of a pointed vertex this rules out
the possibility of having just one face-proper angle.
Indeed, if the face-proper angle is at
the reflex angle then all signs are equal,  and a tangent line
to the point does the job. If the face-proper angle is convex, then a line
through that angle does it. 
Since the number of face-proper angles at a pointed vertex is odd, 
it must be at least three.
\end{proof}

In Theorem~\ref{theoremuniquemax}, we prove an extra condition
that the signs of a self-stress
must satisfy in order to have a non-crossing reciprocal:
the edges in the boundary cycle have opposite sign to those
around the distinguished non-pointed vertex (the one whose reciprocal is
the outer face). The proof uses the relation between self-stresses and
polyhedral liftings of frameworks.


\subsection{Necessary and sufficient conditions}

\begin{figure}[htb]
    \centering
    \includegraphics{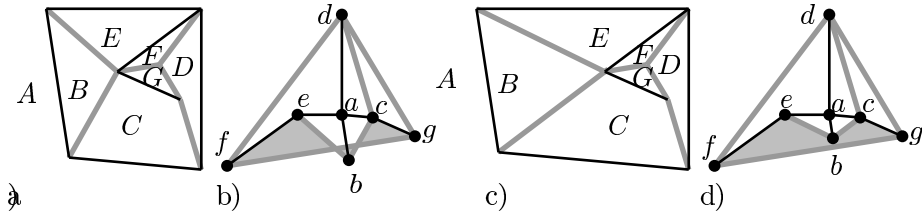}
    \caption{Sign conditions are not enough to guarantee
a planar reciprocal. The geometric circuit in part (a) has the
same signs (and big/small angles) as the one in part (c), but the
former produces a crossing reciprocal while the latter produces a
non-crossing one.\label{figuresignvsframe}}
\end{figure}

Unfortunately, the purely combinatorial conditions on the signs of
the self-stress stated in Theorems~\ref{theorem-face-conditions}
and~\ref{theorem-vertex-conditions}
are not sufficient to guarantee that the reciprocal
is non-crossing. This is illustrated in Figure~\ref{figuresignvsframe}
where a framework and a self-stress satisfying them is shown in (a),
but its reciprocal (b) has crossings.
It has to be noted that this example is a
geometric circuit; in particular it has a unique self-stress (up to a
constant). Part (c) of the same figure
shows a slightly different
embedding of the same graph,  for which the reciprocal
(d) turns out to be
non-crossing. In particular, there can be no purely
combinatorial characterization (that would depend only on the signs of the
self-stress and on which angles are big/small)
of what frameworks have a planar reciprocal.

But it is easy to analyze what goes wrong in this example: the signs
around the interior vertex of degree five in part (a) should produce a
pseudo-quadrangle, 
but instead, it produces a self-intersecting
closed curve (which can be regarded as a ``self-intersecting
pseudo-quadrangle''). The following statement tells us that
such self-intersecting pseudo-quadrangles are
actually the only thing that can prevent a non-crossing framework
with the appropriate signs in its self-stress from having a
non-crossing reciprocal:

\begin{theorem}
\label{badquadrangles}
Let $(G,\rho)$ be non-crossing framework with a given self-stress $\omega$.
The reciprocal  is non-crossing if and only if
the signs of the self-stress
around every vertex satisfy the conditions
of Theorem~\ref{theoremsimultaneous} and, in addition,
the face cycles reciprocal to the non-pointed vertices with
four sign changes are themselves non-crossing \textup(and hence,
pseudo-quadrangles\textup).
\end{theorem}

\begin{proof}
That the reciprocal cycles of every vertex of  $(G,\rho)$
need to be non-crossing for the reciprocal to be non-crossing is
obvious. The reason why we only impose the condition on vertices
of type (2.c) is that the reciprocal cycles of vertices of types (1),
(2.a) and (2.b) are automatically non-crossing: there are
no self-intersecting convex polygons or pseudo-triangles.

Let us see sufficiency. The conditions we now have on vertex cycles
tell us that we have a collection of simple polygons (one of them
exterior to its boundary cycle, containing all the ``infinity'' part
of the plane) and that these polygons can locally be glued to one
another: for every edge of every polygon there is a well-defined
matching edge of
another polygon. Moreover, the fact that the orientations are all
consistent implies that the two matching polygons for a given edge lie
on opposite sides of that edge.

If we glue all these polygons together (which can be done for any reciprocal,
non-crossing or not) what we get is a map from the
topological dual of the framework $(G,\rho)$ (with a point removed,
in the interior of the face reciprocal to the distinguished
non-pointed vertex) to the Euclidean plane. The local argument
shows that this map is a covering map, except perhaps at vertices
where in principle the map could wind-up two or more complete turns.
But the covering map clearly covers infinity once, hence by continuity
it covers everything once. This implies that the map is
actually a homeomorphism. That is, that the reciprocal framework is
non-crossing.
\end{proof}

One can say what the extra condition in Theorem~\ref{badquadrangles}
means for the values of the self-stress in more explicit terms. For
this, observe that a closed cycle is self-intersecting if and only if
it can be decomposed into two cycles. Translated to a vertex
$v$ of $(G,\rho)$ this means that
there are two edges $e$ and $e'$ around $v$ such that
the self-stress (restricted to
the edges around $v$ and the equilibrium condition at $v$)
can be ``split''
into two self-stresses, one supported on the edges on one side of $e$ and $e'$
and another on the edges on the other side. $e$ and $e'$ are allowed
to be used in both self-stresses, but if so with the same sign they
have in the original self-stress.

To see this equivalence, the reader just needs to remember
how to construct the reciprocal cycle  of a given self-stressed vertex:
consider all edges
incident to $v$ oriented going out of $v$ and then place them
one after another (the end of one coinciding with the beginning of
the next one), scaling each edge by the value of the self-stress on that
edge; in particular, reversing the edge if the self-stress is negative.

\medskip

It is interesting to observe that Theorem~\ref{badquadrangles} does not
explicitly require that the framework is a pseudo-quadrangulation,
but instead it gives that as a consequence of the hypotheses.
Corollary~\ref{corovertextoface}
below shows that the vertex conditions alone suffice for this.

\begin{lemma}
\label{lemmavertextoface}
Let $(G,\rho)$ be a pseudo-quadrangulation with $e$ edges, $t$ 
pseudo-triangles,
$q$ pseudo-quadrangles, $y$ pointed vertices and $x$ non-pointed vertices.
In a self-stress of $(G,\rho)$, the following
five properties are equivalent\textup:
\begin{enumerate}
\item the face conditions of Theorem~\ref{theorem-face-conditions}.
\item there are exactly $t$ vertex-proper angles.
\item there are at most $t$ vertex-proper angles.
\item there are exactly $3y+4x-4$ face-proper angles.
\item there are at least $3y+4x-4$ face-proper angles.
\end{enumerate}
\end{lemma}

\begin{proof}
Since the total number of angles $2e$ equals, by 
Lemma~\ref{lemmaquadrangulation},
$t$ plus $3y+4x-4$, conditions (2) and (3) are equivalent to (4) and 
(5), respectively.

For (1)$\Rightarrow$(2)  observe that the face conditions can be rephrased as
``there is exactly one vertex-proper angle in each pseudo-triangle, 
and no vertex-proper
angle in a pseudo-quadrangle or in the outer face''. For the 
converse, (2)$\Rightarrow$(1),  recall that we
always have at least $t$ vertex-proper angles, one at each
pseudo-triangle. The face conditions are just saying that there are no more.
The same observation gives  (2)$\Leftrightarrow$(3).
\end{proof}

\begin{corollary}
\label{corovertextoface}
Let $(G,\rho)$ be a non-crossing framework, and let $\omega$
be a sign assignment satisfying the vertex-condition of 
Theorem~\ref{theorem-vertex-conditions}.
Then, the face conditions of Theorem~\ref{theorem-face-conditions} also hold.
In particular, $(G,\rho)$ is a pseudo-quadrangulation.
\end{corollary}

\begin{proof}
Let $t$ denote the number of pseudo-triangles in $(G,\rho)$, and let 
$q$ be the 
number of other bounded faces
(it will soon follow that they have to be pseudo-quadrangles, but
we don't explicitly require this). The same 
counting argument of
Lemma~\ref{lemmaquadrangulation} yields the inequality 
$2e\ge3t+4q+y$, with equality if and only
if we have a pseudo-quadrangulation. Gluing-in Euler's formula gives 
$2e\le t+3y+4x-4$,
with equality for pseudo-quadrangulations.

Now, the vertex-conditions imply $3y+4x-4$ face-proper angles and we 
have at least $t$ vertex-proper
angles
(one in each pseudo-triangle). Hence, $2e\ge t+3y+4x-4$, which means that
all the inequalities mentioned so far are tight and we have a 
pseudo-quadrangulation
satisfying the face-conditions.
\end{proof}


\subsection{Three special cases}
Something more precise can be said if $(G,\rho)$ is either a geometric
circuit, or a pseudo-triangulation, or if it has a unique non-pointed
vertex. Observe that the first case is self-reciprocal and the other
two are reciprocal to each other.

We start with the case with a single non-pointed vertex.
Recall that the space of self-stresses in a non-crossing framework
is bounded above by the number of non-pointed vertex. Then, one
non-pointed vertex is the minimum needed to sustain a self-stress,
and that self-stress will be unique. We call frameworks with only one
non-pointed vertex \emph{almost pointed}.
The crucial feature about this case is that
the extra condition introduced in Theorem~\ref{badquadrangles}
is superfluous. Actually, in this case the reciprocal is always non-crossing.

\begin{corollary}
\label{corollaryalmostpointed}
Let $(G,\rho)$ be a non-crossing framework with a
single non-pointed vertex and with a self-stress.
Then, the reciprocal framework is non-crossing.
In particular, $(G,\rho)$ is a pseudo-quadrangulation.

Moreover, if the self-stress is everywhere non-zero, then
the reciprocal is a pseudo-triangulation with
$q+1$ non-pointed vertices, where $q$ is the number of pseudo-quadrangles
in $(G,\rho)$.
\end{corollary}

Assuming that the self-stress is everywhere non-zero is actually no
loss of generality: The reciprocal of a singular self-stress is just
the reciprocal of the subgraph of non-zero edges. But we need it in
the second part of the statement in order to get the correct count of 
pseudo-quadrangles.

\begin{proof}
Let $t$ be the number of pseudo-quadrangles and let
$q$ be the number of other faces. As in the previous corollary, we
get that $2e\le t + 3y + 4x-4 = t+3y$, with equality if and only if
we have a pseudo-quadrangulation. But we have equality, since
we have at least $t$ vertex-proper angles (one per pseudo-triangle)
and at least $3y$ face-proper angles, by Lemma~\ref{lemmapointed}.
The equality implies that the vertex conditions are satisfied, hence the
reciprocal is non-crossing.

The fact that the reciprocal is a pseudo-triangulation with $q+1$ 
non-pointed vertices is trivial.
The reciprocals of pointed vertices are pseudo-triangles, and the 
reciprocals of pseudo-quadrangles
and of the outer face are non-pointed vertices.
\end{proof}

\begin{figure}[htb]
   \centering
   \includegraphics{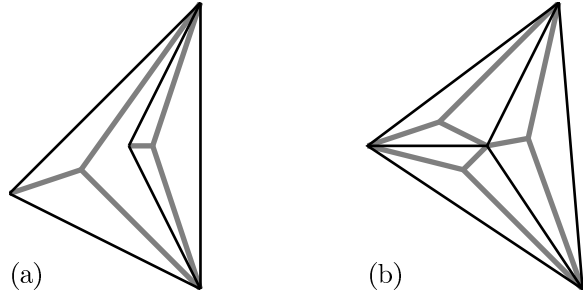}
  \caption{These pseudo-triangulations have a self-stress fulfilling all face
conditions
but do not have a good self-stress.\label{figuretwoconnected}}
\end{figure}

Now we look at pseudo-triangulations. The first observation is that not
all pseudo-triangulations have self-stresses that produce non-crossing
reciprocals. For example, the ones in Figure~\ref{figuretwoconnected}
cannot have self-stresses satisfying the vertex conditions, because  those
conditions forbid more than one non-pointed vertex of degree three.  It is
also interesting to observe that these pseudo-triangulations possess
self-stresses satisfying the face conditions. For example, put negative
stress to all edges incident to non-pointed vertices of degree three, and
positive stress on the others.

\begin{corollary}{\label{corollarypseudo}}  Let $(G,\rho)$
be a pseudo-triangulation with $x$ non-pointed vertices,
and let a self-stress be given to it such that the reciprocal is
non-crossing.
Then, this reciprocal has $1$ non-pointed vertex, $x-1$ pseudo-quadrangles
and $n-x$ pseudo-triangles.
\end{corollary}

\begin{proof}
Straightforward, from
Theorem \ref{theoremsimultaneous}.
\end{proof}

Finally, we look at geometric circuits. Again, not all have 
non-crossing reciprocals, as Figure \ref{figurewheel4recipnonpseudo} 
shows.

\begin{corollary}
\label{corollarygeometric}
Let $(G,\rho)$ be a geometric circuit. That is,
$G$ is a Laman circuit and $(G,\rho)$ is a
non-crossing framework with a non-singular self-stress.

If a reciprocal $G^*$ is non-crossing, then the numbers of
pseudo-triangles, pseudo-quadrangles, pointed vertices and non-pointed
vertices are the same in $G$ and $G^*$.
\end{corollary}

\begin{proof}
We use the formula $2e=t+3y+4x-4$ of Lemma~\ref{lemmaquadrangulation}.
Since a Laman circuit has $e=2 n -2=2x+2y-2$ edges, we get $y=t$.
But $t$ is also the number of pointed vertices in the reciprocal and 
$y$ the number of
pseudo-triangles in it. Now, by Euler's formula, $q+t+n-1=e=2n-2$, hence
$q+t=n-1$ and $q=n-t-1=n-y-1=x-1$. That is to say, $q=x-1$ and $x=q+1$.
Again, $x-1$ is the number of pseudo-quadrangles
in the reciprocal, and $q+1$ the number of non-pointed vertices in it.
\end{proof}


\section{Laman circuit pseudo-triangulations
have planar reciprocals}\label{sectsimultcirc}

\subsection{The non-singular case (geometric circuit pseudo-triangulations)}

Let us start with a geometric circuit pseudo-triangulation. That is 
to say, a Laman circuit
embedded as a pseudo-triangulation with one non-pointed vertex and 
whose self-stress
is non-zero on every edge. This simultaneously satisfies the hypotheses of
Corollaries~\ref{corollaryalmostpointed},~\ref{corollarypseudo} 
and~\ref{corollarygeometric}.
Hence:

\begin{theorem}\label{theoremrecipcircuit}
The reciprocal of a Laman circuit pseudo-triangulation with non-singular
self-stress is non-crossing and again a Laman circuit pseudo-triangulation.
\qed
\end{theorem}

\begin{figure}[htb]
    \centering
    \includegraphics{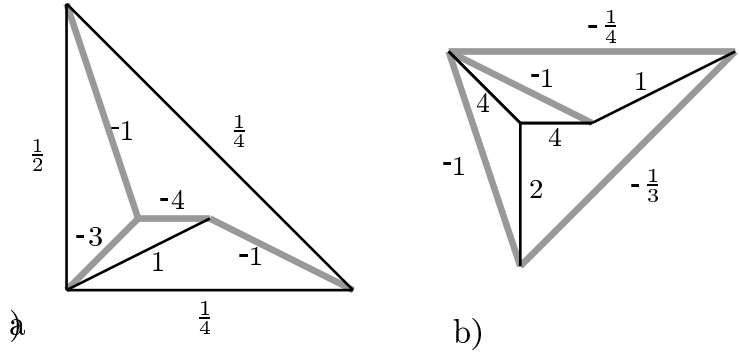}
    \caption{A reciprocal pair of Laman circuit
pseudo-triangulations.\label{figurewheel4recipseudo}}
\end{figure}

Together with Theorem~\ref{theorem2connected},
  Theorem~\ref{theoremrecipcircuit} implies:
\begin{theorem}\label{theoremsimultchararct}
Let $G$ be a Laman circuit. The following are equivalent for $G$:
\begin{enumerate}
\item $G$ is planar.
\item $G$ has a planar embedding with a non-crossing reciprocal.
\item $G$ can be embedded as a pseudo-triangulation with
one non-pointed vertex \textup(whose reciprocal is, in turn,
a pseudo-triangulation with one non-pointed vertex if the embedding is
generic\textup).
\end{enumerate}
\end{theorem}

\subsection{Singular circuit pseudo-triangulations}\label{subsectsing}

For a given planar Laman circuit $G$, the embeddings $\rho$ creating
a pseudo-triangulation $(G,\rho)$ form  an open subset of
$\mathbb{R}^{2|V|}$. The non-singular pseudo-triangulations of
\S\ref{sectsimultcirc}.1, which are a geometric circuit, form an open dense
subset of this subset.  The remaining singular pseudo-triangulations are
`seams' between some components of this open dense set.

\begin{figure}[htb]
   \centering
   \includegraphics{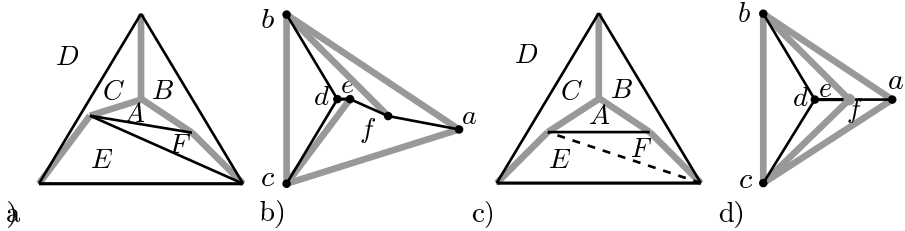}
   \caption{A singular self-stress on a Laman circuit  pseudo-triangulation
can drop an edge from the original graph - and fuse two vertices in the
reciprocal.\label{figurenew1}}
\end{figure}

Consider any singular pseudo-triangulation on a Laman circuit. The self-stress
shrinks to a subgraph $G_s$, see Figure~\ref{figurenew1}a.
Since this framework is still infinitesimally
rigid with $|E|=2|V|-2$, we have a 1-dimensional space of self-stresses.
This framework $(G_s,\rho_s)$  and self-stress can be approached as a
limit of geometric circuits $(G,\rho_n)$ on the whole graph,  each with a
planar reciprocal by \S\ref{sectsimultcirc}.  One can anticipate that the
limit of these reciprocals will also be non-crossing, and this is what we 
prove below.

For example, when  an edge drops out of the self-stress,  the
two faces separated by the lost edge become one face  in the subgraph
$G_s$, see Figures~\ref{figurenew1}(a) and (c).
For each lost edge of the original, the
corresponding edge of the reciprocal, whose length records the coefficient in
the self-stress, will shrink to zero, and the two reciprocal vertices are {\em
fused} into one vertex corresponding to the unified face of the original. See
Figures~\ref{figurenew1}(b) and (d).

We now prove that even in these singular situations the reciprocal
is non-crossing and a pseudo-triangulation.


\begin{theorem}\label{theoremsign} Let $(G,\rho)$ be Laman circuit
embedded as a \textup(possibly singular\textup) 
pseudo-triangulation. Then it has a unique self-stress,
supported on a  subgraph $G_s\subseteq G$. 
 $G_s$ is a pseudo-quadrangulation with a unique non-pointed
vertex and with $q:=2n-2-e$ pseudo-quadrangles,
if $G_s$ has $e$ edges and spans $n$ vertices.\GRcomment{more precise n.}
Its reciprocal is non-crossing, and it is a
pseudo-triangulation with $n-1$ pseudo-triangles and
  $q+1$ non-pointed vertices.
\end{theorem}

\begin{proof}
By Lemma~\ref{lemmaquadrangulation}  $(G,\rho)$ has a unique 
non-pointed vertex.
Clearly, vertices of $G_s$
that were pointed in $(G,\rho)$  are pointed also  in $(G_s,\rho)$. 
In particular, $G_s$ has at most one
non-pointed vertex. By Corollary~\ref{corollaryalmostpointed}, 
the reciprocal
is non-crossing. The other statements are easy.
\end{proof}

\subsection{Good self-stresses.}

We have seen that all reciprocals of a (possibly singular) self-stress
on a Laman circuit are non-crossing pseudo-triangulations.
We say that a self-stress 
on a non-crossing framework $(G,\rho)$ is a \emph{good self-stress} if
it is non-zero on all edges and the reciprocal for this self-stress is
non-crossing.  The existence of a good self-stress is precisely
equivalent to the existence of a non-crossing reciprocal with all
edges of non-zero length.  Does the process of Theorem~\ref{theoremsign}
   create all
examples of a good stress on a pseudo-triangulation?  The answer is
yes.

\begin{theorem}\label{theorempseudorecip}   If a pseudo-triangulation
$(G^*,\rho^*)$ has a good self-stress, then $(G^*,\rho^*)$  is the reciprocal
of a \textup(possibly  singular\textup) Laman circuit 
pseudo-triangulation $(G,\rho)$.
\end{theorem}

\begin{proof}
Let $(G_s,\rho)$ be the reciprocal of $(G^*,\rho^*)$, which is 
non-crossing by assumption.
By Corollary~\ref{corollarypseudo},
$(G_s,\rho_s)$ is an almost pointed framework, with
pseudo-triangles and pseudo-quadrangles.
If $(G_s,\rho_s)$ is already a pseudo-triangulation, then both
graphs are Laman circuits, and we are finished.  Otherwise there
are some pseudo-quadrangles.

As is well-known, a ``diagonal'' edge can be added through the
interior of each pseudo-quadrangle to subdivide it into two
pseudo-triangles, such that no new non-pointed vertex is created
\cite[Theorem 6]{rwwx}.
This process
creates a pseudo-triangulation $(G,\rho_s)$,
with one non-pointed vertex
and the same vertex set as $(G_s,\rho_s)$.

Such an almost-pointed pseudo-triangulation $(G,\rho_s)$  has a
unique self-stress, in this case the self-stress supported on $G_s$.
That is to say, $(G^*,\rho^*)$ is not only the reciprocal of
$(G_s,\rho_s)$, but also of the (singular) self-stress on the
almost-pointed pseudo-triangulation $(G,\rho_s)$.
\end{proof}

\begin{figure}[htb]
   \centering
   \includegraphics[scale=.8]{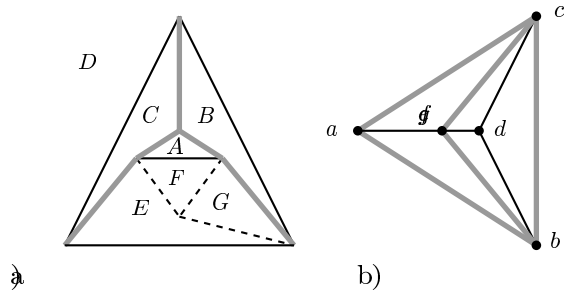}
   \caption{A singular stress on a Laman circuit can drop both vertices
and edges.
\label{figuredropvertex}}
\end{figure}

A singular self-stress can drop not only edges of the original
pseudo-triangulation but also vertices (Figure~\ref{figuredropvertex}).
However,
it is a consequence of this proof that we can choose some alternate
pseudo-triangulation in which the singular self-stress spans all
vertices.  In that case, some simple counting arguments give more
information on the connections between
$(G,\rho_s)$ and the support $(G_s,\rho_s)$ of the singular stress.

\begin{corollary}\label{corollarysign} Let $(G,\rho)$ be a Laman circuit
pseudo-triangulation, and let $G_s$ be a
spanning
subgraph of $G$ supporting
a self-stress.
Let $k$ be the number of edges not used in $G_s$,
so that $|E_s|= 2|V_s|-2 -k$, $k >0$. Then\textup:
\begin{enumerate}

\item $(G_s,\rho)$ is a pseudo-quadrangulation with $n-1-2k$ pseudo-triangles
and $k$ pseudo-quadrangles, each formed as the union of two
pseudo-triangles of $(G,\rho)$.

\item The non-pointed vertex of $(G,\rho)$ is still non-pointed in
$(G_s,\rho)$, and $(G_s,\rho)$ contains the boundary cycle of $(G,\rho)$.

\item The reciprocal is a  pseudo-triangulation with
$k+1$ non-pointed vertices.
\end{enumerate}
\end{corollary}

\begin{proof} The $n-1$ pointed vertices of $G$ are still pointed in $G_s$.
Hence, they still have at least $3n-3$ face-proper angles. Since
every edge is incident to two faces, the removal of $k$ edges destroys at
most $2k$ pseudo-triangles, hence we still have at least $n-1-2k$
pseudo-triangles, each with at least one vertex-proper  angle.
These $(3n-3)+(n-1-2k) = 2(2n-2-k)$ angles equal twice the number of
edges in $G_s$. 
In particular, the number
of pseudo-triangles of $G$ that survived in $G_s$ is exactly $n-1-2k$,
and there is no other pseudo-triangle in $G_s$. Therefore,
each of the $k$ removed edges merged two pseudo-triangles into
a pseudo-quadrangle, and the $2k$ merged pseudo-triangles are all
different. This proves parts (1) and (2) (the latter because if
the removal of an edge makes the non-pointed vertex pointed, then this
removal merges two pseudo-triangles into a pseudo-triangle, not a
pseudo-quadrangle).
\end{proof}

\section{The spatial liftings of non-crossing reciprocal
pairs.\label{sectlift}}

A self-stress on a framework $(G,\rho)$ defines a lifting of it into
3-space  with the property that face cycles are coplanar.  Here we look at
the lifting produced by a \emph{good} self-stress; that is, a self-stress
on a non-crossing framework that produces a non-crossing reciprocal. For
any non-crossing framework, the lifting is a polyhedral surface with exactly
one point above each point of the plane (we have to stress that the outer
face is considered exterior to its boundary cycle;  the standard Maxwell
lifting would consider it interior to the cycle, hence providing a closed,
perhaps self-intersecting surface, with two points above each point
inside the convex hull of the framework). The lifting is unique  up to a
choice of a first plane and of which sign corresponds to a valley or a
ridge
\cite{cw01,cw02,Whiteley2}. Our standard choice for the starting plane places
the exterior face  horizontally at  height zero, and our
standard choice for signs sets the edges in the boundary cycle  as
valleys.  The latter makes sense since in a good self-stress all boundary
edges have  the same sign. We call this the \emph{standard lift} $(G,\mu)$
of the good self-stress
$(G,\rho)$.

Figures~\ref{figurewheel4lpseudolift}, \ref{figurenpr1both},
and~\ref{figurenpr2both} show standard liftings of several good
self-stresses. In all of them one  observes a similar ``shape'':
the entire surface curls 
upwards  from the base to a single
maximum point, which is the lifting of the distinguished non-pointed
vertex.
In particular, there are no local maxima other than
this peak, or local minima except the exterior face in such a surface.
The main theorem in this section shows that these claims hold for all
lifts of non-crossing frameworks with non-crossing reciprocals.

For this spatial analysis
it is easier to reason with the \emph{Maxwell reciprocal}, in which each 
reciprocal
edge is perpendicular, instead of parallel, to the original edge.  This Maxwell
reciprocal is obtained by rotating the Cremona reciprocal by $90^\circ$.
The reason is that 
given a spatial lifting a Maxwell reciprocal is
created by choosing one central point in
$3$-space---for example $(0,0,1)$---and drawing normals to each of the faces
through this point.  The intersection of the normal to face $F$ with the plane
$z=0$ is then the reciprocal vertex for this face \cite{cw01,cw02}. 
If we want
to ensure, for visual clarity, that the reciprocal and the original frameworks
do
not overlap, we simply translate this construction off to the right by picking
the central point to be $(t,0,1)$. This is what we have done in our figures.

  \begin{theorem}\label{theoremuniquemax}
    \begin{enumerate}
    \item [\rm (i)] Given a non-crossing framework $(G,\rho)$ with a
      self-stress such that the corresponding reciprocal $(G^*,\rho^*)$
      is non-crossing, the standard lifting $(G,\mu)$ has a unique
      local \textup(and global\textup) maximum point, whose reciprocal
      is the boundary of $(G^*,\rho^*)$, and has all signs in the
      self-stress opposite to those on the original boundary.

\item [\rm (ii)] The maximum is the unique  point where the lifted surface is
``pointed'', meaning  that a hyperplane exists passing through
it and leaving a neighborhood
of it in one if its two open half-spaces.

    \item [\rm (iii)] The boundary face is the unique local minimum of
      $(G,\mu)$.
    \end{enumerate}
\end{theorem}

\begin{proof} In the standard lifting, the boundary is a horizontal plane and
every edge is attached to an upward sloping face.  No local
maximum can be on the boundary.

Take any isolated local maximum. Cut the lifted surface, just below
this point, with a horizontal plane.  This cuts off a pyramid whose
vertical projection is a non-crossing wheel framework.  Because of the
spatial realization, this projection is a reduced non-crossing
framework $(W,\rho_w)$ (Figure~\ref{figurecutwheels}) which has a
corresponding reduced reciprocal, also, graphically, a wheel.  At this
hub vertex, the pyramid and the original surface have the same face
planes and edges.  Therefore, this hub has the same reciprocal polygon
in the wheel reciprocal and in the original reciprocal $(G^*,\rho^*)$,
and the same stresses along these edges in the two projections.
Moreover, the signs of the spokes in this stress indicated the
concavity or convexity of the edge in the lifting $\mu$ and the
concavity or convexity in the rim polygon at this spoke of the wheel.
We can now do a sweep around the maximal hub: we start with the plane of
a face, then we rotate the plane about an adjacent spoke until we reach the
next face, etc.  The normals to these tangent planes track the
reciprocal polygon, with a reciprocal vertex for the normal to each of
the faces.
\begin{figure}[htb]
   \centering
   \includegraphics{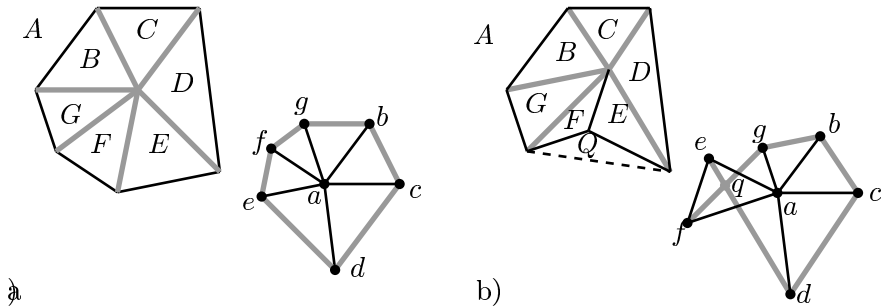}
   \caption{Possible pyramids which might be created by a
cut near a maximum in the lifting of a non-crossing
framework (shown in vertical projection), with their reciprocals.}
\label{figurecutwheels}
\end{figure}

If the base of the pyramid is not a convex polygon, then there is some
segment of the convex hull between two vertices of the base polygon
which does not lie in the polygon, placing all other vertices into one
half-plane (Figure~\ref{figurecutwheels}b).

Consider the plane $Q$ formed by this segment and the hub in space.
This plane will place all of the pyramid in one half-space, touching
it in at least two spokes.  We claim that the normal $q$ to this plane
is a crossing point of the reciprocal polygon.  As we sweep around a
spoke that lies in~$Q$, we will encounter $Q$ as one of the tangent
planes, and hence $q$ will appear on the edge reciprocal to the spoke.
Since this happens for at least two spokes, the reciprocal polygon is
self-intersecting at $q$.

Since we know that the reciprocal polygon is non-crossing (and non-touching) we
conclude that the original pyramid base must be a convex polygon.  A direct
analysis of
the wheel now guarantees that all the signs of the spokes of the wheel are the
same---and are opposite to the rim of the wheel, representing ridges. 
Since the
stress, and the signs, at the hub are the same in the larger framework, we
conclude that any maximal vertex has all signs the same, and these signs are
opposite to the boundary of the framework.

By Theorem~\ref{theoremsimultaneous}, we know that there is only one 
vertex with
no
sign changes in each side of a reciprocal pair of non-crossing  frameworks, and
that this vertex is the reciprocal of the boundary of the other 
framework in the
pair. We conclude that there is a unique local maximum---the global maximum.
This global maximum corresponds to the convex boundary polygon in the
reciprocal.
Thus we have proved part~(i) of the theorem for a local
maximum consisting of an isolated vertex.

If we can cut off a vertex $v$ of the lifting with a plane $P$ which is
not necessarily horizontal, but cuts all edges between $v$ and its
neighbors, the above argument applies without change. 
This section will still produce a spatial wheel which will
project into a plane wheel with the hub inside the rim. Any such vertex
would have to be the reciprocal vertex of the boundary, and since this
vertex is unique, we have proved statement~(ii) of the theorem.

Let us consider the case that a local maximum is not just a vertex but
a larger connected set $M$ consisting of horizontal lines and horizontal faces.
We can apply the argument of the previous paragraph to
  any vertex of the convex hull of~$M$.

Given any local minimum that is not on the boundary, the same argument about
the pyramid cut, and the sense of the reciprocal polygon applies.  This gives a
contradiction, so the only  local minimum is along the boundary face,
proving part~(iii).
\end{proof}

Statement~(ii) can be interpreted as saying that the maximum is the
only (locally) strictly convex vertex of the surface.  However, it is
possible for a vertex $v$ to have all adjacent vertices in a
\emph{closed} half-space through $v$.  This can only happen for a
pointed vertex and the half-space bounded by the plane of the face
into which it points will for a pointed vertex. In particular, all
boundary vertices have this property.

We can derive the following consequences of the previous theorem:

\begin{corollary}
\label{corollaryopposite}
In a good self-stress, the edges incident to the distinguished 
non-pointed vertex have sign
opposite to that of the boundary cycle.
\end{corollary}

\begin{proof}
Since the distinguished vertex is the maximum of the standard 
lifting, some (hence all) of its adjacent
edges are ridges in it. In the contrary, the boundary edges are 
valleys by definition
of the standard lifting.
\end{proof}

\begin{corollary}
\label{corollarylevelcurves}
In the lifting of a  good self-stress there are no 
(horizontal) saddle points. All 
level curves are simple closed curves, and as we increase the height 
the level curve moves monotonically from the boundary cycle to the 
distinguished non-pointed vertex.
\end{corollary}

The reader has to observe, however, that the intermediate level 
curves need not be convex,
even tough the boundary cycle and the level curves  sufficiently 
close to the tip are convex.
This happens, for example, in Figure~\ref{figurenpr1both}.

\begin{proof}
There cannot be any saddle points, because general Morse theory
on a disc with a horizontal boundary shows that a saddle point
would require an additional local maximum or a new local minimum.
In the absence of saddle points, Morse theory also implies that all level
curves are isotopic to one another, hence they are all simple closed curves
because the boundary cycle is.
\end{proof}

\begin{proposition}
\label{corollary-twisted-saddle}
 For any vertex except the maximum, there is a  plane through that
  vertex that cuts the neighborhood in the ``saddle-point'' way into $4$
  pieces.  
No plane through a
  vertex cuts the neighborhood into more than $4$
  pieces, i.~e., there are no ``multiple saddles''.  

More precisely, for every general direction in the interior of the reciprocal
figure there is a unique vertex such that a plane with this normal through
the vertex cuts the neighborhood in the ``saddle-point'' way into $4$
pieces.  For planes with this normal, all other vertex neighborhoods
are cut into two pieces, with the exception of the peak, whose neighborhood
lies entirely below the plane. For all directions in the outer face
of the reciprocal, a plane with that normal will cut the neighborhood
of every vertex, including the peak, into two pieces.
\end{proposition}

\begin{proof}
The neighborhood of $v$ consists of an alternating cyclic sequence of
edges emanating from $v$ and faces between those edges; in the face
$V^*$ reciprocal to $v$, these correspond to vertices and edges,
respectively.

For a given plane $Q$ through $v$, we want to count how many faces
incident to $v$ it intersects.
In order to determine whether it intersects a given face $F$ between
two neighboring edges $e_1$ and $e_2$ emanating from $v$, we look at
the relation between the reciprocal normal vector $q^*$ and the
reciprocal vertex $f^*$, forming an angle of $V^*$ with the incident
edges $e_1^*$ and $e_2^*$.
It turns out that $Q$ intersects $F$ if and only if
\begin{itemize}
\item 
 $e_1^*$ and $e_2^*$ have different signs and
 the line through $q^*$ and $f^*$  ``cuts through'' the boundary
 of $V^*$ at $f^*$,  or 
\item
 $e_1^*$ and $e_2^*$ have the same sign and
 the line through $q^*$ and $f^*$ is ``tangent'' to the boundary
of $V^*$ at $f^*$.
\end{itemize}
We know that the signs around $V^*$ have to satisfy the face
conditions. It follows
that $q^*$ has the
above-mentioned relation to precisely 4 vertices of a face $V^*$ if
$q^*$ lies in $V^*$, and to precisely 2 vertices of $V^*$ if $q^*$
lies outside $V^*$.
This holds for the interior faces of the reciprocal (pseudo-triangles
and pseudo-quadrangles), and it can be easily proved by checking a few
elementary cases and then showing that the number of ``related''
vertices does not changes as one moves $q^*$, except when crossing the
boundary of~$V^*$.

When $v$ is the peak and $V^*$ is the outer face, $q^*$ is related
to precisely 2 vertices of $V^*$ if it 
lies in the outer face, and to no vertices otherwise.

The desired statements now follow easily.
(The first part of the Theorem, which is only a local statement about
the neighborhood of $v$, can also be proved directly in the original
framework along the lines of the proof of
Theorem~\ref{theoremuniquemax} by considering the geometry and the
possible sign patterns of edges between two intersections with~$Q$.)
\end{proof}

The behavior exhibited in the previous corollary is analogous to what
happens in the upper half of the \emph{pseudosphere}, which is a 
surface in 3-space which serves as a
model for (a part of) the hyperbolic plane, and has constant negative
curvature everywhere.
{The pseudosphere is the surface of revolution generated by
a tractrix. The upper half is given in parametric form
by the equations
$z=u-\tanh u$, $r=\mathop{\rm sech} u$ in polar coordinates 
$r,\varphi,z$, for $u\ge 0$.}
When the pseudosphere is viewed as the graph of a function
over the unit circle, then for every
gradient direction there is a unique point with that gradient, and
the mapping  point  $\leftrightarrow$ gradient  is an
orientation-reversing
mapping between the pseudosphere and the plane.
Even the properties in Theorem~\ref{theoremuniquemax} 
confirm that the lifted surface of each framework in
a non-crossing pair has  the  shape of a rough piecewise-linear
pseudosphere,
except that the pseudosphere has the vertical axis as an asymptote,
whereas our lifted surface reaches a finite maximum.
(In this sense the surface 
$z=(r-1)^2$ over the region $r\le 1$ might be a more appropriate
smooth model for our lifted surface. It has a constant negative
Laplacian $\Delta z=-2$.)
  In a visible sense, this lifted surface is as non-convex as possible.


Liftings of pseudo-triangulations have also
emerged recently in the context of \emph{locally convex} piecewise linear
functions over a
polygonal domain subject to certain height restrictions~\cite{aabk}.
There are some similarities, in particular regarding
the twisted saddle property discussed in
Proposition~\ref{corollary-twisted-saddle},
but we have not explored whether
there are any deeper connections to our present work.

\section{Open Problems}\label{sectopen}

Throughout this section we say that a non-crossing framework is {\em good}
if it has some good self-stress, i.~e., an everywhere
non-zero stress that produces a non-crossing reciprocal.
In particular, every almost pointed non-crossing framework is good.
We do not define a framework to be good if \emph{all} its
self-stresses are good, simply because every framework with at least two
(linearly independent) self-stresses has some bad
self-stresses. Indeed, let $\omega_1$ and $\omega_2$ be two
independent everywhere non-zero
self-stresses on a framework, and suppose they are both
good. Consider the associated standard liftings,
$\mu_1$ and $\mu_2$. For $c>0$ sufficiently big (resp., sufficiently
small) the lifting $c\mu_1$ lies completely above (resp.,
completely  below) $\mu_2$. Since $c\mu_1$ is never equal to $\mu_2$,
there must be an intermediate value $c_0$ for which $c_0\mu_1$ has
some parts above and some parts below $\mu_2$. In particular, 
$c_0\omega_1-\omega_2$ cannot be a good self-stress, because its
associated
lifting has parts above and below the plane of the outer face (which
contradicts Theorem~\ref{theoremuniquemax}).

\subsection{Good and bad pseudo-triangulations}
\begin{figure}[htb]
   \centering
   \includegraphics[scale=.8]{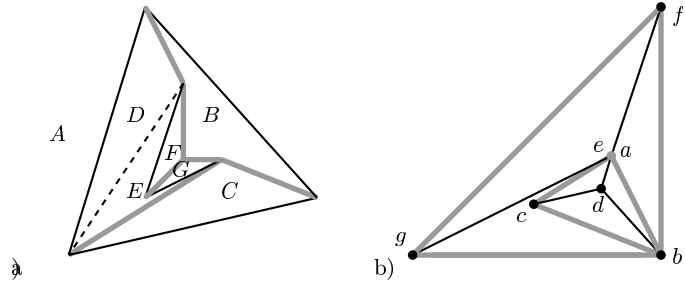}
   \caption{A  pseudo-triangulation with two
non-pointed vertices (b) with its
  reciprocal (a),  a singular
pseudo-triangulation.\label{figuretwonew2}}
\end{figure}

In \S\ref{sectsimultaneously} we have seen examples of pseudo-triangulations
which cannot hold a good self-stress.
In Figure~\ref{figuretwonew2}(b) we see a different framework on
the same graph
                    as in Figure~\ref{figuretwoconnected},
with the
same face structure, which does have a good self-stress,
as demonstrated by the non-crossing reciprocal.
The difference between these examples lies in the choice of
exterior face and of the associated big and small angles.
In~\cite{osss} a combinatorial analog of pseudo-triangulation
is discussed in which the planar graph has a formal labelling of 
`big' and `small'
assigned to the angles, which may or may not correspond to the big 
and small angles of
any planar realization, with pseudo-triangle and pointedness describe 
combinatorially.

For a given graph, having a good self-stress
is an open property, so if there are any good realizations
there are nearby generic realizations.

\begin{open}
What conditions on a combinatorial
pseudo-triangulation ensure existence of  a good generic
realization?

If there is a good generic realization of
a combinatorial pseudo-triangulation,
are all generic realizations of this combinatorial
pseudo-triangulation good? If the answer is yes, then a purely
combinatorial characterization of 
combinatorial pseudo-triangulations that admit 
good embeddings should exist. Find it.
\FScomment{added the last part. GR. last sentence sounds a bit martial}
\end{open}

Observe that the second question has a negative answer if posed
for pseudo-quadrangulations.
Figures~\ref{figurenpr1both} and \ref{figurepseudoquadbad} show
two different generic embeddings of a Laman circuit
as pseudo-quadrangulations with the same big and small angles. The unique
self-stress is good in the first embedding and bad in the second.

\subsection{Good pseudo-quadrangulations.}
Our characterization in  \S~\ref{sectsimultaneously} applies to all
non-crossing reciprocal pairs based on a graph $G$ and its dual $G^*$ 
on the induced
face structure of the embedding.   In general, these reciprocal pairs are
composed of pseudo-triangles and pseudo-quadrangles, with at least one
non-pointed vertex.

\begin{figure}[htb]
   \centering
   \includegraphics{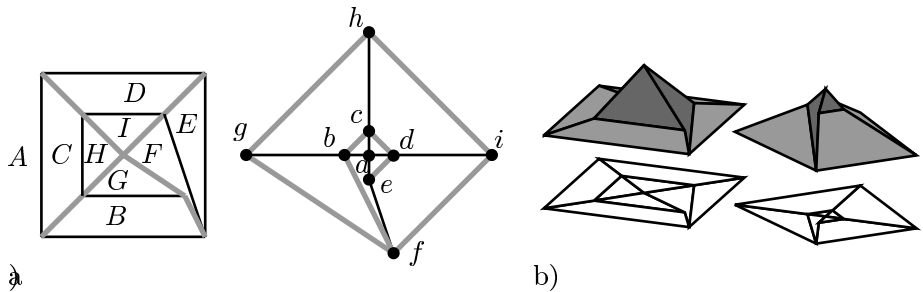}
   \caption{A reciprocal pair of pseudo-quadrangulations on a Laman
circuit.\label{figurenpr1both}}
\end{figure}

\begin{figure}[htb]
   \centering
   \includegraphics{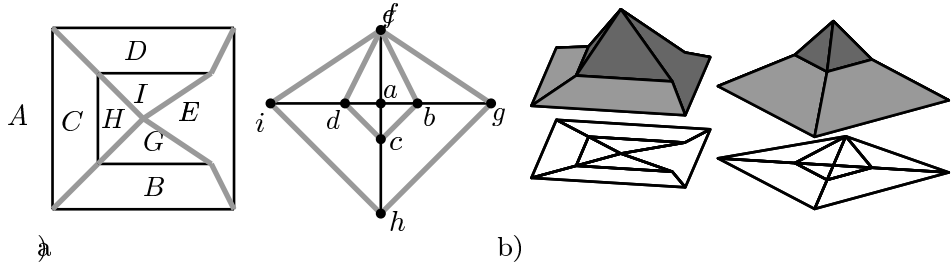}
   \caption{A nearby singular framework with its
compressed reciprocal and lift.\label{figurenpr2both}}
\end{figure}

\begin{figure}[htb]
   \centering
   \includegraphics{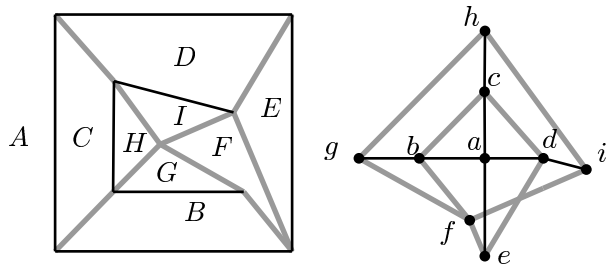}
   \caption{The same graph and combinatorial pseudo-quadran\-gu\-lation can be bad,
because of an altered sign pattern. \label{figurepseudoquadbad}}
\end{figure}

Figure~\ref{figurenpr1both}(a)
contains such a pair of non-crossing reciprocals on a Laman circuit, 
with the sign
pattern of the corresponding self-stress, as predicted by
Theorem~\ref{theoremsimultaneous}, and Figure~\ref{figurenpr1both}b
shows the lifts guaranteed by Maxwell's Theorem.
Figure~\ref{figurenpr2both}(a)
contains a singular framework and its  reciprocal (with two vertices
fused) found as a limit from the previous example, with the sign pattern
of the corresponding self-stress, as predicted by Theorem 3, and
Figure~\ref{figurenpr2both}(b)
shows the lifts.  We note that a small additional change in location, 
away from the
original pair across this singularity, can make this 
pseudo-quadrangulation bad as the
sign pattern is altered on the singular edge, see
Figure~\ref{figurepseudoquadbad}.

This illustrates that, once we leave the realm of
pseudo-triangulations, the big and small angles are not sufficient
to determine whether even a generic framework on a Laman circuit 
has a good self stress. 
\WWcomment{Modified this paragraph}

\begin{open}
\label{open1} Characterize directly by their geometric
properties as embeddings,
all  non-crossing frameworks on a  Laman
circuit
$G$ with a non-crossing reciprocal.
\end{open}

Like all characteristics of a self-stress, the existence of a good 
self-stress on a
framework is invariant under any \emph{external projective transformation
of the plane}, that is,  a projective transformation in which the 
line being sent
to infinity does not intersect any of the vertices or edges of the
original framework. This  invariance is a direct consequence of the map on the
stress-coefficients induced by such a projective transformation, which does not
change any signs and preserves all equilibria \cite{rothwhiteley}. These two
properties of the map guarantee that all the necessary and sufficient 
conditions
of \S\ref{sectsimultaneously} are preserved.

However, whatever projective
property is required for a good self-stress, it is more refined than the simple
choice of large and small angles, or even the oriented matroid of the vertices
themselves.  We do not have a firm conjecture for the geometry in the framework
which determines that there is a good self-stress.  Such a good self-stress
will be singular on any refinement of the framework to a pseudo-triangulation.
The pure condition polynomials mentioned in the introduction, which 
are zero for
edges dropped in a singular self-stress, are projective invariants
\cite{WhiteWhiteley2}. Recall Figure~\ref{figurepseudoNonCircuit}, where the
condition  was the concurrence of three lines.  It is possible that 
all the necessary
information lies in the geometry of singular self-stresses on a
pseudo-triangulation refining the pseudo-quadrangulation.



\subsection{Refinement}
It is natural to start with a non-crossing Laman circuit
pseudo-triangulation and generate denser examples by adding edges. But
the refinement process does not preserve ``goodness'' in
general, even if it is made with no addition of vertices.

\begin{figure}[htb]
   \centering
   \includegraphics[scale=.8]{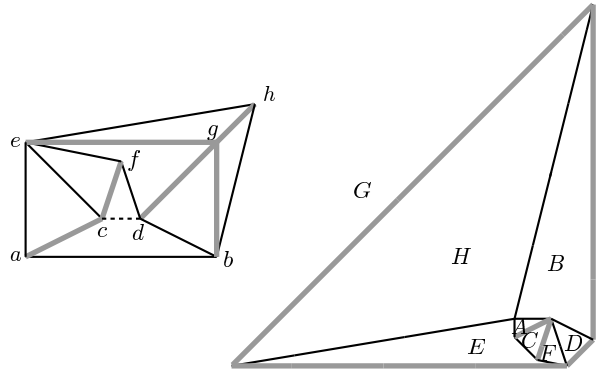}
   \caption{A pseudo-quadrangulation with a good self-stress for which
the added edge (dashed) makes a good self-stress on all edges
impossible.\label{figurecounter}}
\end{figure}

In Figure~\ref{figurecounter} we show an example of this.
The framework on the left, without the edge $cd$, is an
almost-pointed pseudo-triangulation, whose
unique self-stress is by our results good (its reciprocal is on the
right).  We claim that if we 
add the edge $cd$ then there is no good everywhere-non-zero
self-stress on the framework. To see this, let $\omega_1$ be the
self-stress of the pseudo-triangulation and let $\omega_2$ be the
self-stress with support on the subgraph induced by the vertices $a$,
$b$, $c$, $d$, $e$ and $g$ (including the edge $cd$). Since this
subgraph can be lifted to a roof-like surface, $\omega_2$ has one
sign, say negative, on the boundary and the opposite sign in the
interior. W.l.o.g. we assume that $\omega_1$ was also negative on the
boundary. It turns out that $\omega_1$ and $\omega_2$ generate the
space of self-stresses of the whole graph, but no linear combination
$\alpha\omega_1 + \beta\omega_2$ will give alternating signs to the
edges
of the quadrilateral $abdc$: 
if $\alpha/\beta>0$ then $bd$ and $cd$
get
the same sign; 
if $\alpha/\beta<0$ then $ac$ and $cd$ do. 

However, if a good pseudo-triangulation $(G,\rho)$ is refined to 
another pseudo-triangulation $(G',\rho)$ without adding vertices,
we conjecture that $(G',\rho)$ is good too. Our reason for this
conjecture is that it is easy to prove that, at least, there is a
self-stress in $(G',\rho)$ satisfying the vertex conditions of
Theorem~\ref{theorem-vertex-conditions}. The idea of the proof is that
one can go from $(G,\rho)$ to $(G',\rho)$ via a sequence of
``elementary refinements'' of one of the following types:
addition of a single edge
to divide a pseudo-triangle into two, or addition of  
three edges forming a triangle that separates the three corners of a
pseudo-triangle, dividing it in four pseudo-triangles.
(See \cite[Lemma~3]{rwwx} for a formal argument that these refinements
are sufficient.) In both cases, a simple case study together with
elementary properties of self-stresses gives the required sign pattern.

Our conjecture reduces only to prove 
that the extra condition ruling out bad quadrangles can 
be obtained too.

\begin{open}
    Is it true that every pseudo-triangulation refinement (with no
extra vertices) of a good pseudo-triangulation is good?
\end{open}

And if one allows something more than simply adding edges:
\begin{open}
    To what extent can non-crossing reciprocal pairs be generated from Laman
    circuit pseudo-triangulations via refinement?
\end{open}

Solving any of these two problems could be a step towards the first
problem we posed.


%

\subsection{Lifting Questions}
 Liftings of pseudo-triangulations
have also emerged recently in the context of ``locally convex''
functions over a polygonal domain subject to certain height
restrictions~\cite{aabk}.
We have
seen here that, given a pair of frameworks whose reciprocals
are both non-crossing, the lifted surface of each
framework shares some characteristic properties of a pseudo-sphere. It would
be of interest to see if this resemblance increases with the
density of the framework.
\begin{open}
     Let $\{(G_i,\rho_i)\}$ be a sequence of Laman circuit
pseudo-trian\-gu\-la\-tions such that each
     $(G_{i+1},\rho_{i+1})$ is obtained from  $(G_i,\rho_i)$ by
a making a Henneberg II move in
     a randomly selected face.  With an appropriate
normalization of the stresses, this defines a sequence of
lifted surfaces, all of which
     have negative discrete curvature at every vertex except for
a single maximum vertex.

Does this process converge to some limit?
What can be said about the limiting surface? Are there
combinatorial conditions on the
sequence of frameworks that ensure that the limit is something like a
smooth pseudo-sphere?
\end{open}

One could ask the same question with a different model of generating
``random'' Laman circuits.  For example, the PPT-polytope
of~\cite{rss} is a polytope whose vertices are in one-to-one
correspondence with the pointed pseudo-triangulations on a given point
set.  Choosing an extreme vertex in a random direction, (for a
randomly generated point set) produced a pseudo-triangulation, to which
we can add an edge to create a ``random'' Laman circuit.
Limit shapes like this have appeared in other contexts, for example
for random convex polygons~\cite{b,brsz}.

The results in \S\ref{sectlift} interpreted
necessary conditions in \S\ref{sectsimultaneously}  on the self-stress
at the vertices and faces of the original non-crossing framework as
necessary conditions on the lifting.
It is natural to reverse this idea by starting from a piecewise-linear
surface and projecting it back to the plane. Which properties
of a surface are necessary and sufficient?

\begin{open}
\label{open2} Characterize geometrically exactly  which
piecewise-linear spatial surfaces, projecting 1-1 onto the
entire plane, project to planar frameworks with non-crossing
reciprocals.
\end{open}

\subsection{Non-crossing reciprocals with respect to other embeddings}
If a graph  admits several topologically different embeddings
in the plane, one may decide to  construct
  the reciprocal of a non-crossing framework on that graph taking as 
face and vertex cycles
those of a different plane embedding from the one given by the framework.
Or one may not allow this but decide to call the reciprocal 
non-crossing if it is non-crossing as a
geometric graph, even if its face structure is not dual to the one in 
the original.

Our characterizations  of non-crossing reciprocal pairs do not 
address these situations.
Note that this is only an issue for non-3-connected graphs, and that
graphs with cut vertices need not be considered: every self-stress 
can be decomposed
as a sum of self-stresses supported in 2-connected components.

There are two questions that should be addressed here.
\begin{open}
Is there a non-crossing framework $G$ (necessarily $2$-connected but 
not $3$-connected)
whose natural reciprocal
has no crossing edges but is embedded differently from the graph-theoretic
planar dual of $G$?
\end{open}
We know that this cannot happen when $G$ is a Laman circuit (see
Theorem~\ref{theoremsimultchararct}),
but for frameworks with more edges, the question is open.

\begin{open}
   Characterize pairs of non-crossing frameworks which are reciprocals 
to one another,
but not necessarily with respect to the face and vertex cycles given by their
embeddings as frameworks.
\end{open}

\subsection{What planar graphs produce non-crossing reciprocal pairs?}
We finish with perhaps the broadest question of all:

\begin{open}
Given a $2$-rigid planar graph, decide (give a characterization, or at
least a reasonable  algorithm) whether  there is a non-crossing generic
embedding of it that has a good self-stress.
\end{open}

We know for example that for all Laman circuits such an embedding 
exists: any generic
pseudo-triangulation embedding works. In order for a generic 
framework to have a a self-stress on all edges,
it must be $2$-rigid---remain rigid after deletion of any one
edge \cite{JordanJackson}.  But we also know that not all
planar $2$-rigid graphs have such an embedding
(Figure~\ref{figuretwoconnected}).

\section{Acknowledgements}
This research was initiated at the Workshop on Rigidity Theory and
Scene Analysis organized by Ileana Streinu at the Bellairs
Research Institute of McGill University in Barbados, Jan.\ 11--18,
2002 and partially supported by NSF grant CCR-0203224.

\end{document}